\documentclass[]{IEEEtran}
\usepackage{amsmath,amsfonts}
\usepackage{algorithmic}
\usepackage{algorithm}
\usepackage{array}
\usepackage[caption=false,font=normalsize,labelfont=sf,textfont=sf]{subfig}
\usepackage{textcomp}
\usepackage{stfloats}
\usepackage{url}
\usepackage{verbatim}
\usepackage{graphicx}
\usepackage{cite}
\usepackage{xcolor}

\usepackage[short]{optidef}
\usepackage{amsthm}
\newtheorem*{assumption*}{\assumptionnumber}
\providecommand{\assumptionnumber}{}
\makeatletter

 \newtheorem{assume}{Assumption}
 \newtheorem{remark}{Remark}

\makeatletter

\begin{document}

\title{A Comparative Study of Distributed Feedback Optimizing Control Architectures}

\author{Risvan Dirza$^1$, Hari Prasad Varadarajan$^2$, Vegard Aas$^1$, Sigurd Skogestad$^1$, and Dinesh Krishnamoorthy$^{2,*}$

\thanks{$^1$Risvan Dirza, Vegard Aas, and Sigurd Skogestad are with the Department of Chemical Engineering, Norwegian University of Science and Technology (NTNU), NO-7491 Trondheim, Norway (e-mail: \texttt{rdir@equinor.com, vgard.aas@hotmail.com; sigurd.skogestad@ntnu.no}).}
\thanks{$^2$Dinesh Krishnamoorthy and Hari Prasad Varadarajan are with the Department of Mechanical Engineering, Eindhoven University of Technology, 5600 MB, Eindhoven, The Netherlands (e-mail: \texttt{d.krishnamoorthy@tue.nl, h.p.varadarajan@tue.nl}).}
\thanks{R.D and S.S. acknowledges funding support from SFI SUBPRO under Grant Number 237893, which is financed by the Research Council of Norway, major industry partner and NTNU. D.K. acknowledges funding support from the Dutch Research Council (NWO) through the VENI talent scheme (Project number: 20172). }
\thanks{$^*$ Corresponding author.}}

\markboth{Submitted to IEEE Transactions on Control Systems Technology}%
{Dirza \MakeLowercase{\textit{et al.}}: A Comparative Study of Distributed Feedback-optimizing Control Schemes: An Experimental Validation}

\maketitle

\begin{abstract}
This paper considers the problem of steady-state real-time optimization (RTO) of interconnected systems with a common constraint that couples several units, for example, a shared resource. Such problems are often studied under the context of distributed optimization, where decisions are made locally in each subsystem, and are coordinated to optimize the overall  performance. Here, we use distributed feedback-optimizing control framework, where the local systems and the coordinator problems are converted into feedback control problems. This is a powerful scheme that allows us to design feedback control loops, and estimate parameters locally, as well as provide local fast response, allowing different closed-loop time constants for each local subsystem. This paper provides a comparative study of different  distributed feedback optimizing control architectures using two case studies. The first case study considers the problem of demand response in a residential energy hub powered by a common renewable energy source, and compares the different feedback optimizing control approaches using simulations. The second case study experimentally validates and compares the different approaches using a lab-scale experimental rig that emulates a subsea oil production network, where the common resource is the gas lift that must be optimally allocated among the wells. 
\end{abstract}


\section{Introduction}
\IEEEPARstart{F}{eedback}-optimizing control, a class of optimization methods originally introduced in \cite{Morari1980StudiesStructures}, is a powerful technique that transforms the economic optimization problem into a feedback control problem. This method has garnered significant research attention and widespread application across various fields, including process control and autonomous systems \cite{Krishnamoorthy2022Real-TimeReview,hauswirth2024optimization}. The key reason for its popularity lies in its ability to maintain (near-)optimal system performance by promptly (on fast timescales) responding to changing conditions and uncertainties using simple feedback tools such as Proportional Integral Derivative (PID) controllers and selectors.

Often the optimal operation corresponds to case when some of the constraints are active. In such cases, optimal operation can be achieved by active constraint control, where some of the available degrees of freedom can be used to control the constraint to its limiting value using feedback controllers. The unconstrained degrees of freedom and the changing set of active constraints requires special treatment, whereby one has to select controlled variables, known as \textit{self-optimizing} variables, which when kept at constant setpoint, leads to (near-)optimal operation at steady-state \cite{Skogestad2000PlantwideStructure}. The ideal self-optimizing variable is in fact the steady-state cost gradient with respect to the unconstrained degrees of freedom, which when controlled to a constant setpoint of zero, satisfies the necessary condition of optimality \cite{Halvorsen2003OptimalVariables}. The idea of estimating the steady-state cost gradient and driving it to zero using feedback control is also the backbone of several different approaches, including extremum seeking control \cite{Ariyur2003Real-TimeControl}, NCO tracking control \cite{Srinivasan2008TrackingFunction}, Neighboring extremal control \cite{gros2009optimizing} to name a few. The main difference in these approaches lies in how the steady-state cost gradient is estimated. See for example  \cite{Krishnamoorthy2022Real-TimeReview} for a comprehensive overview of model-based and model-free steady-state gradient estimation schemes that can be used in the context of feedback optimizing control.

The works mentioned above primarily concentrated on bridging the optimization problem and the control problem by identifying and estimating the optimal self-optimizing controlled variables. However, when dealing with real-world processes, many of them are large-scale in nature with several interconnected subsystems, posing practical challenges in constructing feedback-optimizing control structures at such a scale. Individual user preferences,   autonomy, and privacy also necessitate the need for local decision-making within each unit without the need to share local information in the form of models, real-time measurements, operational schedules etc. 
Hence, it becomes imperative to integrate decomposition methods \cite{Lasdon1970OptimizationSystems,Chiang2007LayeringArchitectures} into the framework of feedback-optimizing control. In this context, feedback optimizing controllers are designed locally for each susbsytem, which are coordinated to achieve  optimal performance of the overall system.  We refer to this class of problems as \emph{distributed feedback-optimizing control}, wherein both the control and estimation techniques can be decomposed and designed independently within each subsystem.

In the distributed optimization literature, there are two main
strategies that can be used to decompose a large-scale problem
into several smaller subproblems, namely,
into: 
\begin{itemize}
\item Primal decomposition - which is an opportunity cost-based  framework. In this framework, the master coordinator directly allocates the shared resource to subsystems, which report their incurred opportunity costs (also known as marginal costs) associated with the allocated resource. Subsequently, the allocation of the shared resource is adjusted to ensure equal opportunity cost across all subsystems.
\item Dual decomposition - which is a price-based coordination framework. In a price-based coordination framework, subsystems interact through a market-like mechanism where prices of the shared resources are adjusted by a central coordinator to achieve market equilibrium. Subsystems then adjust their resource usage based on these prices.
\end{itemize}
 In general, the distributed feedback-optimizing control scheme holds several crucial benefits for interconnected process control systems (in order of the importance):
\begin{enumerate}
	\item This scheme enables formulating each local problem independently (including independent local estimator), assuming that the overall objective consists of local objectives that are additively separable.
	\item Fast local responses are achievable because one may use different closed-loop time constants for each local system.
	\item Solutions can be implemented immediately after each iteration, eliminating the need to wait for convergence over several iterations.
	\item This scheme circumvents slow response and numerical robustness issues like divergence, often encountered in Real-Time Optimization (RTO) methods that solve a  numerical optimization problem online.
	\item The distributed feedback-optimizing control scheme serves as an alternative to numerical-based RTO approaches, which may be computationally expensive.
\end{enumerate}
Recently, we proposed  distributed feedback optimizing control approaches based on the two decomposition frameworks listed above to design  control architectures to achieve optimal steady-state operation of the overall system using only feedback  controllers in each  subsystem \cite{Krishnamoorthy2021ASharing,Dirza2021OptimalOptimization,Dirza2022ExperimentalRig,Dirza2022Real-timeDecomposition,dirza2024primal}. 
While these methods have been studied and developed independently, a comprehensive comparative study analyzing the advantages and disadvantages of the different distributed feedback optimizing control architectures is notably lacking.

\paragraph*{Main contribution} This paper builds upon our previous works \cite{Krishnamoorthy2021ASharing,Dirza2021OptimalOptimization,Dirza2022ExperimentalRig,Dirza2022Real-timeDecomposition,dirza2024primal} to conduct a comparative analysis of three distributed feedback optimizing control architectures. We leverage two case studies, one experimental and one simulation-based, to provide a comparative study of the three distributed feedback optimizing control architectures. Our key contributions are as follows.
\begin{itemize}
\item Our previous work \cite{Dirza2022ExperimentalRig} experimentally validated the feedback optimizing control structure based on dual decomposition on a lab-scale gas-lifted oil well rig. One of the contributions of this paper is the experimental validation of the control  architectures resulting from primal decomposition and  dual decomposition with override, which has not been done before.
\item We then compare the three schemes of distributed feedback-optimizing control using the experimental test facility.
\item This paper also  considers a  different case study  that aims to achieve \textit{demand response} within a residential energy hub powered by renewable energy to further demonstrate and compare  the  different distributed feedback optimizing control architectures. The use of distributed  feedback optimizing control for demand response management has not been studied before. \end{itemize}
The reminder of the paper is organized as follows. Section \ref{Sec: Problem Formulation} lays the groundwork by describing the three distributed feedback optimizing control architectures in detail.
Section~\ref{sec:Demand response} present a comparative study of the three control architectures using simulations on a demand response case study.  Section~\ref{Sec: Lab Setup} presents a comparative study of the three control architectures using experiments conducted in a lab-scale test facility that emulates a subsea oil production network. Section~\ref{sec: Discussion} provides a comprehensive comparison of the three control architectures across various aspects, drawing insights from both the simulation and experimental results, before concluding the paper in Section~\ref{Sec: Conclusion}.


\section{Problem Formulation for Distributed Feedback-optimizing Control}\label{Sec: Problem Formulation}
Consider an integrated steady-state optimization problem consisting  of $N$  subsystems, coupled by a shared resource constraint
\begin{mini!}
	{\{\mathbf{u}_{i}\}}{ \sum_{i=1}^N J_i\left(\mathbf{u}_i,\mathbf{d}_i\right)}
	{\label{Eq:gen_prob}}{}\label{Eq:gen_prob_obj}
	\addConstraint{\sum_{i=1}^N\mathbf{g}_i\left(\mathbf{u}_i,\mathbf{d}_i\right) \leq\bar{\mathbf{g}}}
	\label{Eq:gen_prob_cons}
\end{mini!}
where for each subsystem, $\mathbf{u}_{i}\in \mathbb{R}^{n_i}$ are the set of decision variables, $\mathbf{d}_{i}  \in  \mathcal{D}_i \subset \mathbb{R}^{d_{i}}$ denotes the set of disturbances with $\mathcal{D}_i$ being a compact set. $J_{i}:\mathbb{R}^{n_i} \times \mathbb{R}^{d_{i}} \rightarrow \mathbb{R}$ is the local cost function, and $\mathbf{g}_{i}:\mathbb{R}^{n_i}\times \mathbb{R}^{d_{i}} \rightarrow \mathbb{R}^{m}$ denotes the (coupling) constraints with  $\bar{\mathbf{g}} \in \mathbb{R}^{m}$ being upper limit of the $m$ shared resources that must be optimally allocated among the $N$ subsystems. Further we  assume that the functions $J_i$ and $\mathbf{g}_i$ are continuous in both $\mathbf{u}_{i}$ and $\mathbf{d}_i$ and $\exists \; \{\mathbf{u}_i\}$ such that \eqref{Eq:gen_prob_cons} is feasible for all $\{\mathbf{d}_{i}  \in  \mathcal{D}_i\}$. Such problems appear in a broad range of engineering applications, such as energy management systems, chemical process plants, water distribution networks, advanced manufacturing systems to name a few. 

\begin{remark}
We aim to find the optimal value for the control variables $\mathbf{u}$ for a given  disturbance $\mathbf{d}$ that will drive the system to a desired steady-state operating point. At steady-state, the rate of change of the states $\dot{\mathbf{x}}$  becomes zero, essentially eliminating the need for the state variables themselves in the optimization problem \eqref{Eq:gen_prob}. While states are not explicitly included, the cost  and constraint functions  are often implicitly dependent on the states. 
\end{remark}

This paper considers the following problem setting:
\begin{enumerate}
	\item Each subsystem $i =1,\dots,N$  optimizes its own objective $J_{i}(\mathbf{u}_{i},\mathbf{d}_{i})$ locally without sharing its local models and real-time measurements with the other subsystems, while coordinating to achieve the optimum of \eqref{Eq:gen_prob}.
	\item Each subsystem uses only simple feedback controllers such as PI controllers, and does not solve numerical optimization problems online.
\end{enumerate}
To satisfy the first requirement, we use distributed optimization \cite{Lasdon1970OptimizationSystems} to decompose the large-scale optimization problem \eqref{Eq:gen_prob} into smaller subproblems, that can be coordinated to achieve overall optimum. 

To satisfy the second requirement of achieving optimal steady-state operation without solving numerical optimization problems online, we consider feedback optimizing control, where the optimization objectives are translated to feedback control objectives. For example, this would typically involve controlling the gradient of the cost function to a constant setpoint of zero, thereby satisfying the necessary condition of optimality. In the context of interconnected subsystems, optimizing the individual subsystems (by controlling the gradient of their respective local cost functions $\nabla_{\mathbf{u}_i} J_i$ to zero) does not aggregate to optimal operation of the overall system. In this case, one has to control a self-optimizing variable, that also depends on the coordinating information that can be used to coordinate
the local feedback controllers, thus achieving system-wide optimal operation. 
Doing so would result in a distributed feedback optimizing control structure, where each subsystem achieves the necessary conditions of optimality for a given coordinating information, and as the  coordinator updates the coordinating information, we achieve optimal steady-state operation of the the overall system.  There are different approaches to designing the control architectures, which are described in the subsequent subsections.
\begin{figure*}
	\centering
	\includegraphics[width=\linewidth]{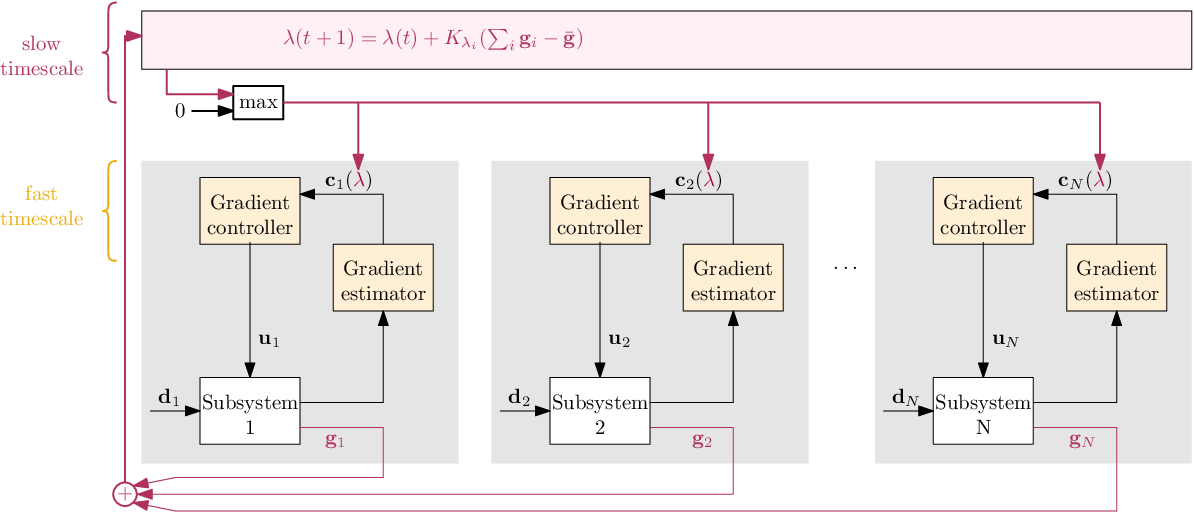}
	\caption{Schematic representation of the distributed feedback optimizing control architecture based on dual decomposition (price-based coordination).  Gray box indicates the information boundary of each subsystem.}\label{Fig:Dual_CSD}
\end{figure*}
\paragraph*{Notational remark} The collection of decision variables and the disturbances of the overall system is denoted by  $\mathbf{u}:= \{ \mathbf{u}_{i}\}_{i=1}^N$ and $\mathbf{d}:= \{ \mathbf{d}_{i}\}_{i=1}^N$, respectively. For the sake of notational simplicity $\{\cdot\}$ by default indicates $\{\cdot\}_{i=1}^N$, and $\|\cdot\|$ by default indicates the Euclidean norm, unless specified otherwise. The element-wise product of two vectors of equal dimension $ \mathbf{a},\mathbf{b} \in \mathbb{R}^m$ is denoted by the Hadamard product $ \mathbf{a}\circ \mathbf{b}$.


\subsection{Control structure architecture using price-based coordination}
In price-based coordination, the individual subsystems adjust their local decision-making  in response to varying price signals to collectively achieve the supply-demand equilibrium in this micro-market setting. The  actions of the individual subsystems, guided by market prices, leads to a coordinated outcome where resources are allocated efficiently.

The basic idea here is to relax the coupling constraint using Lagrangian relaxation, and solve the resulting dual problem:
\begin{equation}\label{Eq:duality}
	\max_{\lambda\ge0} \; \min_{\mathbf{u}}\; 	\mathcal{L}(\mathbf{u},\mathbf{d},{\lambda})
\end{equation}
where
\begin{equation}\label{Eq:Lagrangian}
	\mathcal{L}(\mathbf{u},\mathbf{d},{\lambda}) = \sum_{i=1}^{N} J_{i}(\mathbf{u}_{i},\mathbf{d}_{i}) + \lambda^{\mathsf{T}}\left( \sum_{i=1}^{N}\mathbf{g}_{i}(\mathbf{u}_{i},\mathbf{d}_{i}) - \bar{\mathbf{g}}\right)
\end{equation}
with $\lambda \in \mathbb{R}^m_{\geq0}$ being the Lagrange multiplier for the shared resource constraint \eqref{Eq:gen_prob_cons}. 



For a given value $\lambda(t)$ at time $t$, the stationarity condition of \eqref{Eq:duality} is the solution of the following unconstrained problem,
\begin{mini}
	{\mathbf{u}}{\mathcal{L}\left(\mathbf{u},\mathbf{d},\lambda(t)\right)}{q(\lambda(t)):=}
	{\label{cen_prob_primal}}{}
\end{mini}
which is expressed as,
\begin{equation}\label{Eq:Stationarity}
    \nabla_{\mathbf{u}} \mathcal{L}(\mathbf{u},\mathbf{d},\lambda(t)) = 0
\end{equation}
Therefore, to achieve optimal operation at steady-state, we need to control $\mathbf{c}(\lambda(t))\ := \nabla_{\mathbf{u}} \mathcal{L}(\mathbf{u},\mathbf{d},\lambda(t))$ to a setpoint of zero.
Since the cost and constraints are additively separable, the controlled variable
\begin{equation}\label{Eq: dual_decomp}
	\begin{aligned}
		\mathbf{c}(\lambda(t))\ &:= \nabla_{\mathbf{u}} \mathcal{L}(\mathbf{u},\mathbf{d},\lambda(t))\\ &= \begin{bmatrix}
			\nabla_{\mathbf{u}_1} J_1(\mathbf{u},\mathbf{d})\\\vdots\\\nabla_{\mathbf{u}_N} J_N(\mathbf{u},\mathbf{d})
		\end{bmatrix} + \lambda(t)^\top \begin{bmatrix} \nabla_{\mathbf{u}_1} \mathbf{g}_{1}(\mathbf{u},\mathbf{d}) \\ \vdots \\ \nabla_{\mathbf{u}_N} \mathbf{g}_{N}(\mathbf{u},\mathbf{d}) \end{bmatrix}
	\end{aligned}
\end{equation}
can be easily decomposed, where each subproblem now controls
\begin{equation}\label{Eq:ci}
	\mathbf{c}_i(\lambda(t))\ :=  \nabla_{\mathbf{u}_i} J_i(\mathbf{u},\mathbf{d}) + \lambda(t)^\top \nabla_{\mathbf{u}_i} \mathbf{g}(\mathbf{u},\mathbf{d})
\end{equation}
to setpoint of $\mathbf{c}_i^{sp}=0$ by using $\mathbf{u}_i$ as the manipulated variables in the inner loop. Here $\mathbf{c}_{i} \in \mathbb{R}^{n_{i}}$ is the gradient of the local subproblem.

The coordinator problem is then given by $\max_{\lambda\ge0} q(\lambda)$, which is a maximization problem with $m$ degree of freedom. The price of the shared resource (i.e., Lagrange multiplier) $\lambda(t)$ is updated by the central coordinator  using a dual ascent step. That is, at every sample time of the central coordinator $ t,t+1,\dots $, the total resource consumption $\sum_{i=1}^{N}\mathbf{g}_{i}$ is reported to the central coordinator, which is used to update the shadow price as follows.
\begin{equation}\label{Eq: max_dual}
	{\lambda}(t+1) = \text{max}\left[0,{\lambda}(t) + K_{{\lambda}_i} \left(\sum_{i=1}^{N}\mathbf{g}_{i}(\mathbf{u}_{i},\mathbf{d}_{i}) - \bar{\mathbf{g}} \right)\right]
\end{equation}
Note that the price update step \eqref{Eq: max_dual} is simply an integral controller  with a  \emph{max} selector that controls the coupling constraint \eqref{Eq:gen_prob_cons} in a slower timescale in a cascaded fashion.
\begin{figure*}
	\centering
	\includegraphics[width=\linewidth]{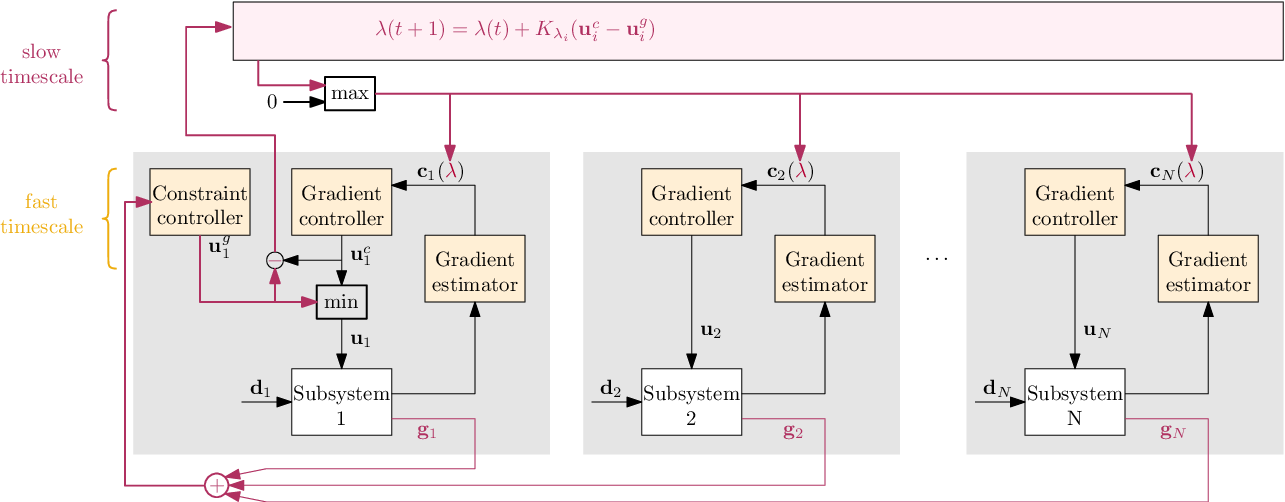}
	\caption{Schematic representation of the distributed feedback optimizing control architecture based on dual decomposition (price-based coordination) with active constraint control override in the fast timescale.  Gray box indicates the information boundary of each subsystem.}\label{Fig:DualOverride_CSD}
\end{figure*}

While controlling \eqref{Eq:ci} to a setpoint of 0 in each subsystem ensures the satisfaction of the stationarity condition at steady-state, the integral controller with the \emph{max} selector in \eqref{Eq: max_dual} satisfies the remaining KKT conditions, namely
\begin{subequations}
    \begin{align}
    \sum_{i=1}^{N}\mathbf{g}_{i}(\mathbf{u}_{i},\mathbf{d}_{i}) - \bar{\mathbf{g}}& \leq 0\\
    \lambda & \ge 0 \\
    \lambda \circ  \left( \sum_{i=1}^{N}\mathbf{g}_{i}(\mathbf{u}_{i},\mathbf{d}_{i}) - \bar{\mathbf{g}}\right) &= 0
\end{align}
\end{subequations}


\begin{remark}\label{rem:ADMM}
	One can also use an augmented Lagrangian function with quadratic penalty terms instead of the plain Lagrangian function in \eqref{Eq:Lagrangian} to guarantee  convergence of the coordinated control structure. See  \cite{Krishnamoorthy2021ASharing} for  the detailed convergence analysis.
\end{remark}
To summarize, the control structure based on price-based coordination uses local feedback controllers in each unit  to control the variable $ \mathbf{c}_{i}(\lambda)  $ to a constant setpoint of  $\mathbf{c}_i^{sp}=0$, for a given price $ \lambda $, and in the layer above, a central coordinator updates the price of the shared resource $\lambda$ to match the demand and the supply of the shared resource in the slow timescale. 



\begin{remark}\label{rem:slow constraint control}
	When using price-based coordination, the shared resource  constraint is controlled in the slow time scale using the shadow prices $\lambda$. Since this is done in the slow time scale in the central coordinator, whose closed-loop time constant is at least 5 times slower than the closed-loop time constants of the local controllers, the shared resource constraint will be satisfied only in the slow time scale. In general, controlling the constraint in the slow timescale may not be desirable, especially if these are hard constraints, or if the cost of constraint violation is high.
\end{remark}

\subsection{Control structure architecture using price-based coordination with constraint  override}

In general, it is desirable to control the constraints in the faster timescale, and look after optimality in the slower time scale (i.e., first feasibility, then optimality). However as remarked above, price-based coordination scheme accounts for feasibility in the slow timescale, which may not be desirable. To address this challenge, we recently proposed to move the constraint control to the faster timescale by overriding some of the local controllers with the constraint control whenever a coupling constraint exceeds its limit \cite{Dirza2022Primal-dualControl,dirza2024primal}.

Here, for all the subsystem $i=1,\dots,N$,  local controllers are designed to control the gradient $\mathbf{c}_i({\lambda}(t))$ to a setpoint of $\mathbf{c}_i^{sp}={0}$ using $\mathbf{u}_{i}$ as the manipulated variable.  Let us denote the control input computed by the gradient controller as $\mathbf{u}_{i}^c$.
In addition, for a pre-selected subset of local subsystems, which we call, \emph{critical subsystems},  an overriding constraint control is designed that controls the constraint $ \sum_{i=1}^{N}\mathbf{g}_{i}(\mathbf{u}_{i},\mathbf{d}_{i}) - \bar{\mathbf{g}}$ to a setpoint of $0$ using $\mathbf{u}_{i}$ as the manipulated variable (active constraint control). Let us denote the control input computed by the constraint controller as $\mathbf{u}_{i}^g$. Since two controlled variables are paired with the same manipulated variable $\mathbf{u}_{i}$ for the critical subsystems, a  selector block is used to switch between the two controllers as described in  \cite{Krishnamoorthy2020SystematicSelectors}. 
For example, a minimum selector block is used $\mathbf{u}_i = \min(\mathbf{u}_{i}^c,\mathbf{u}_{i}^g)$ to reduce the consumption of the shared resource and thus ensure constraint feasibility in the fast timescale.

If the coupling constraint is violated, the gradient controllers in the critical subsystems are overridden by the active constraint controller, and in the slow timescale, the shadow prices are updated such that the gradient controller is selected by the selector block $\mathbf{u}_i = \min(\mathbf{u}_{i}^c,\mathbf{u}_{i}^g)$ at steady-state (i.e., we want $\mathbf{u}_{i}^c\leq\mathbf{u}_{i}^g$). This can be obtained by designing an integral control in the coordinator that controls the \textit{auxiliary constraint} $(\mathbf{u}_{i}^c-\mathbf{u}_{i}^g) \leq 0$ by updating the shadow prices as shown below.
\begin{equation}\label{Eq: max_override}
	{\lambda}(t+1) = \text{max}\left[0,{\lambda}(t) + K_{{\lambda}}(\mathbf{u}_{i}^c-\mathbf{u}_{i}^g)\right]
\end{equation}
A \emph{max} selector is used to ensure non-negative Lagrange multipliers.
\begin{figure*}
	\centering
	\includegraphics[width=\linewidth]{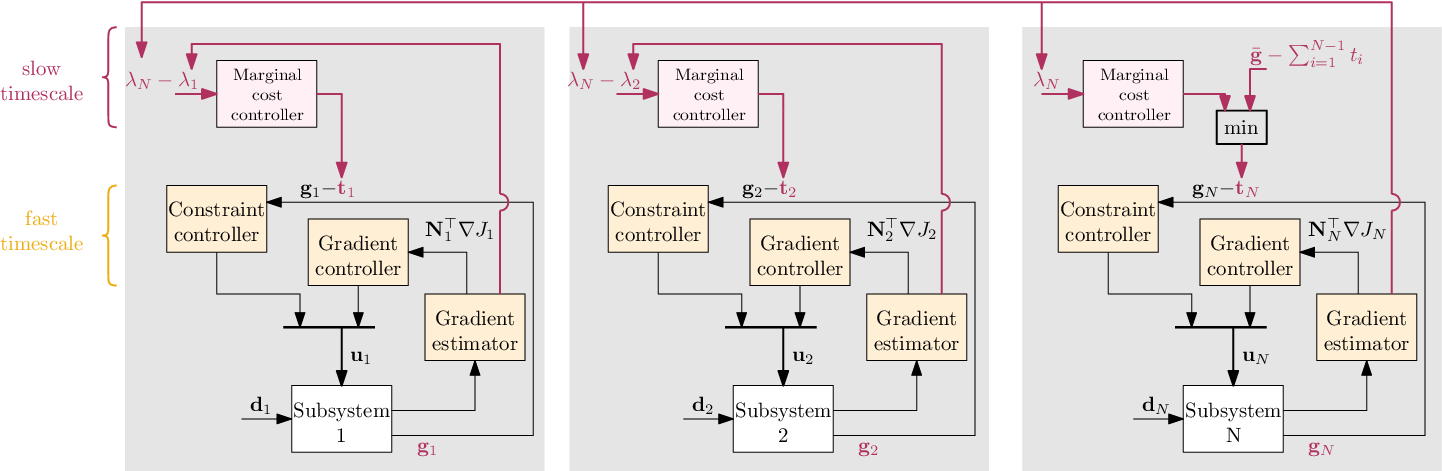}
	\caption{Schematic representation of the distributed feedback optimizing control architecture based on primal decomposition (opportunity cost based coordination). Grey box indicates the information boundary of each subsystem. }\label{Fig:Primal_CSD}
\end{figure*}
To summarize, this provides a price-based coordination framework, where the coupling constraints are controlled in the fast timescale. Although controlling the coupling constraints in the fast timescale using  override can help minimize  dynamic constraint violation, it cannot ensure hard constraint feasibility.
Furthermore, moving the constraint control to the fast timescale requires one to make a pairing decision for the constraint controllers in the fast timescale.

\subsection{Control structure architecture using  opportunity cost based coordination}

A straightforward approach to ensure that the coupling constraint remains feasible also during the transients, is to directly allocate the available shared resource to the different subsystems. Here, the shared resource is directly allocated in the slow timescale based on the local opportunity cost reported by the different subsystems, such that the opportunity cost of all the units are the same. 
Such a decomposition and coordination framework is  known as primal decomposition or decomposition by right hand side allocation \cite{Bertsekas2016-az,Lasdon1970OptimizationSystems}.

Here the local problem for each subsystem is given by
\begin{subequations}\label{Eq:Primal subroblem}
	\begin{align}
		V_i(\mathbf{t}_{i}) = 	\min_{\mathbf{u}_{i}} &\; J_{i}(\mathbf{u}_{i},\mathbf{d}_{i}) \\
		\textup{s.t.} &\; \mathbf{g}_{i}(\mathbf{u}_{i},\mathbf{d}_{i}) = \mathbf{t}_{i}  \label{Eq:primal constraint}
	\end{align}
\end{subequations}
where the $ \mathbf{t}_{i}  \in \mathbb{R}^m$ is the amount of the shared resource allocated to the $ i^{th} $ subsystem.

\begin{assume}\label{asm:sufficientMVs}
For each subsystem there are at least as many manipulated variables as the constraints, i.e.,  $n_{i} \geq m$ for all $i=1\dots,N$.
\end{assume}
To transform  subproblems \eqref{Eq:Primal subroblem} to a feedback control problem,
\begin{enumerate}
	\item Using the first $m$ degrees of freedom, control the local constraint $\mathbf{g}_{i}(\mathbf{u}_{i},\mathbf{d}_{i}) - \mathbf{t}_{i} $ to a setpoint of 0.
	\item Using the remaining $n_{i}-m$ unconstrained degrees of freedom, control the reduced gradient $ \mathbf{N}_i^\top \nabla_{\mathbf{u}_i} J_i\left(\mathbf{u}_i,\mathbf{d}_i\right) \in \mathbb{R}^{n_{i}-m}$, to setpoint of $\mathbf{0}$, where $\mathbf{N}_i^\top$ is the null space of $\nabla_{\mathbf{u}_i} \mathbf{g}_i\left(\mathbf{u}_i,\mathbf{d}_i\right)$ \cite{Jaschke2012OptimalSystems,Krishnamoorthy2020LinearVariables}.
\end{enumerate}
To this end, the local subproblems control \[ \mathbf{c}_{i}(\mathbf{t}_{i}) = \begin{bmatrix}
	\mathbf{g}_{i}(\mathbf{u}_{i},\mathbf{d}_{i}) - \mathbf{t}_{i} \\
	 \mathbf{N}_i^\top \nabla_{\mathbf{u}_i} J_i\left(\mathbf{u}_{i},\mathbf{d}_{i}\right)
\end{bmatrix} \] to a setpoint of $\mathbf{c}_{i}^{sp} = 0$ for a given allocation of $\mathbf{t}_{i}$ in the fast timescale.
The cost incurred by fixing the constraint to the pre-allocated share of the common resource $ 	\mathbf{g}_{i}(\mathbf{u}_{i},\mathbf{d}_{i}) - \mathbf{t}_{i}  $ is the \emph{opportunity cost}, which  represents  the potential benefit or value that the subsystem  could have gained by relaxing the constraint. In other words, this can be interpreted as the implicit prices the different units are willing to pay for additional resources.  This provides insight into how changes in resource availability affect the overall system's objective function, which can be used by the central coordinator to optimally allocate the shared resource.  The opportunity cost, also known as the marginal cost, is simply the Lagrange multiplier $ \lambda_i $ corresponding to the constraint \eqref{Eq:primal constraint}. Using the necessary condition of optimality of the local subproblems $ \eqref{Eq:Primal subroblem} $, this can be computed by solving the system of linear equations
\begin{equation}\label{Eq: local lang}
\nabla_{\mathbf{u}_i}J_i\left(\mathbf{u}_i^*,\mathbf{d}_i\right) + \nabla_{\mathbf{u}_i}\mathbf{g}_i\left(\mathbf{u}_i^*,\mathbf{d}_i\right) \lambda_i^* = 0
\end{equation}

If $m< n_i$, the constraint gradient $\nabla_{\mathbf{u}_i}\mathbf{g}_i\left(\mathbf{u}_i,\mathbf{d}_i\right)$ is not a square matrix. Denoting the cost and constraint gradients as $\Gamma_i:= \nabla_{\mathbf{u}_i}J_i\left(\mathbf{u}_i^*,\mathbf{d}_i\right) $ and $\Phi_{i}:= \nabla_{\mathbf{u}_i}\mathbf{g}_i\left(\mathbf{u}_i^*,\mathbf{d}_i\right) $, respectively,  we can use the pseudo-inverse to compute the opportunity cost as:
\begin{equation}\label{Eq:OpportunityCost}
	 \lambda_i^* = -  \left[\Phi_{i}^{\mathsf{T}}\Phi_{i}\right]^{-1}\Phi_i^{\mathsf{T}}\Gamma_i, \quad \forall i = 1,\dots,N.
\end{equation}
\begin{remark}\label{rem:Square}
	In the special case of $n_{i}=m$, the reduced gradient controller is not required, i.e., $\mathbf{c}_{i}(\mathbf{t}_{i}) =
		\mathbf{g}_{i}(\mathbf{u}_{i},\mathbf{d}_{i}) - \mathbf{t}_{i} $.
\end{remark}
To ensure that the total consumption  is less than the supply at all times, the resource allocated to the last unit is given by
\begin{equation}\label{Eq:Primal Feasible}
	\mathbf{t}_{N} = \bar{\mathbf{g}}- \sum_{i=1}^{N-1} \mathbf{t}_{i}
	\end{equation} Simply put, the remaining shared resource is allocated to subsystem $N$. 
	The coordination  problem then involves finding the optimal resource allocation $ \mathbf{t}_{i} $ for $ i = 1,\dots,N-1 $, which can be expressed as
	\begin{align}\label{Eq:Primal master}
		\min_{\{\mathbf{t}_{i}\}_{i=1}^{N-1}} & \sum_{i=1}^NV_i(\mathbf{t}_{i})
	\end{align}
 where $V_i(\mathbf{t}_{i}) := J_i\left(\mathbf{u}_i^*,\mathbf{d}_i\right) +  \lambda_i^{*^\mathsf{T}}\left(\mathbf{g}_i\left(\mathbf{u}_i^*,\mathbf{d}_i\right) - \mathbf{t}_i \right)$.
	Given \eqref{Eq:Primal Feasible}, the necessary condition of optimality of \eqref{Eq:Primal master} is satisfied when $ \lambda_N -\lambda_i = 0$ for all $ i= 1,\dots,N-1 $, i.e., equal marginal cost for all the units. Therefore in each unit $ i= 1,\dots,N-1 $, we use feedback controllers  in the slower timescale to control the local marginal cost $\lambda_N-\lambda_i  $ to a setpoint of 0 using $\mathbf{t}_{i}$ as the manipulated variable. For  $i=N$, $\mathbf{t}_N$ is given by \eqref{Eq:Primal Feasible}. 
 
 If the resource constraint is not optimally active (i.e., more supply than demand for the shared resource), the allocated resource will be sufficient for all the subsystems, and consequently  the opportunity costs for all the subsystems will be equal to zero. To achieve this automatically, we  use a feedback controller on the $ N^{th} $ subsystem that controls $ \lambda_N$ to a setpoint of 0 using $\mathbf{t}_{N}$ as the manipulated variable, and use a minimum selector block to choose the minimum value between the controller output and \eqref{Eq:Primal Feasible} as shown in Fig.~\ref{Fig:Primal_CSD}.
	Since the other $ N-1 $ units control $ \lambda_N -\lambda_i$ to zero, controlling $ \lambda_N $ to zero, automatically results in  $ \lambda_i =0$ for all $ i = 1,\dots,N $, which from \eqref{Eq: local lang} implies that $\nabla_{\mathbf{u}_i}J_i\left(\mathbf{u}_i,\mathbf{d}_i\right)=0$ for all $i$.

\begin{remark}
	The problem can be equivalently decomposed using primal decomposition by introducing a slack $\mathbf{g}_{0} \geq 0$ to transform the inequality coupling constraint to an equality coupling constraint, such that $\mathbf{g}_{0} + \sum_{i=1}^{N}\mathbf{g}_{i} - \bar{\mathbf{g}} = 0$. This can be interpreted as a virtual subsystem, whose  corresponding opportunity cost is $\lambda_0 = 0$, since subsystem $i=0$ does not contribute towards the overall cost. In this case, controlling $\lambda_N - \lambda_0$ to a setpoint of 0 using $\mathbf{t}_{0}$, and setting  $\mathbf{g}_{0} = \max(0,\mathbf{t}_{0})$ is equivalent to controlling the marginal cost $\lambda_N$ to a setpoint of 0 in the $N^{th}$ subsystem.
\end{remark}

\section{Case Study 1: Demand response in a district-level energy hub}\label{sec:Demand response}
In recent years, the integration of renewable energy sources, such as solar power, into district-level energy systems has gained significant traction as a means to promote sustainability and reduce reliance on traditional fossil fuel-based energy generation. Demand response (DR) has emerged as a pivotal strategy to enhance the flexibility and resiliency of  energy hubs integrated with renewable energy sources.
For example, in the context of temperature control and electric heating, demand response would involve adjusting the settings on the heating system such as lowering the desired temperature slightly during times when there's strain on the  grid.
\begin{figure}
	\centering
	\includegraphics[width=\linewidth]{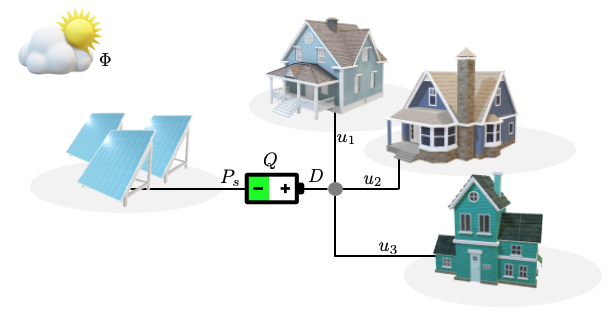}
	\caption{Schematic representation of a district-level micro energy hub  powered by a common solar energy source coupled with battery storage}\label{Fig:EnergyHub}
\end{figure}
This case study investigates a district-level energy hub designed for residential sustainability. The hub leverages a centralized solar power source coupled with battery storage to provide renewable energy for electric heating systems in participating houses. This configuration empowers residents to utilize clean energy while storing surplus energy for later use, fostering a more sustainable community energy ecosystem.

Integrating solar and battery storage at the district level presents unique considerations for demand response (DR) initiatives within the energy hub. The intermittent nature of solar generation and the limitations of battery capacity necessitate dynamic management of energy consumption to guarantee reliable heating for residents. However, the flexibility of battery storage opens doors for optimizing energy use and responding to real-time fluctuations in supply and demand.

Demand response programs aim to achieve a more efficient balance between electricity supply and demand, optimizing grid operation, minimizing strain during peak periods, and facilitating the integration of renewable energy sources. These programs typically involve modifying electricity consumption patterns, such as shifting usage to off-peak hours or reducing overall consumption during peak demand or grid constraints. In this case study, we explore the implementation of demand response strategies using various distributed feedback optimizing control architectures.
\begin{figure}
	\centering
	\includegraphics[width=0.99\linewidth]{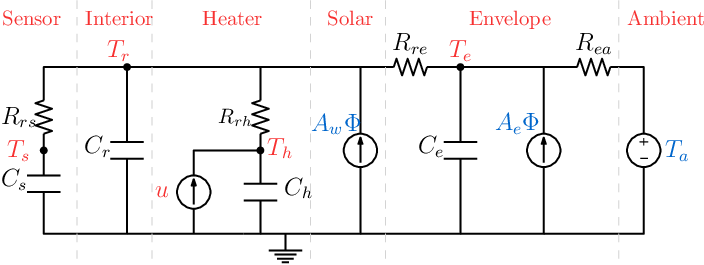}
	\caption{RC-thermal model of a single house.}\label{Fig:HouseRC}
\end{figure}
We consider a small-scale district where $N=3$ residential buildings share power from a common solar energy source with battery storage to meet their heating needs, as shown in Fig.~\ref{Fig:EnergyHub}. Each house is modelled as a RC-thermal model \cite{bacher2011identifying} (represented in the form of a circuit diagram in Fig.~\ref{Fig:HouseRC}).
	\begin{align*}
	\frac{\textup{d} T_{s}}{\textup{d}t} & = \frac{1}{R_{rs}C_{s}}(T_{r} - T_{s})\\
	\frac{\textup{d} T_{r}}{\textup{d}t} & = \frac{1}{R_{rs}C_{r}}(T_{s} - T_{r}) +\frac{1}{R_{rh}C_{r}} (T_{h}-T_{r}) + \frac{A_{w}\Phi}{C_{r}}  \\
	& \qquad +\frac{1}{R_{re}C_{r}} (T_{e}-T_{r}) + \frac{1}{R_{ra}C_{r}} (T_{a}-T_{r})\\
	\frac{\textup{d} T_{h}}{\textup{d}t} & = \frac{1}{R_{rh}C_{h}}(T_{r} - T_{h}) + \frac{u}{C_{h}}\\
	\frac{\textup{d} T_{e}}{\textup{d}t} & = \frac{1}{R_{re}C_{e}}(T_{r} - T_{e}) + \frac{1}{R_{ea}C_{e}}(T_{a} - T_{e}) + \frac{A_{e}\Phi}{C_{e}}
\end{align*}
The subscripts $ (\cdot)_s,(\cdot)_r,(\cdot)_h, (\cdot)_e $ and $ (\cdot)_{a} $ represent the sensor, building interior, heater, building envelope, and ambient, respectively. Here, $ T $ denotes temperature in [$^\circ$C], $ R $ denotes thermal resistance in [$^\circ$C/kW] , $ C $ denotes heat capacity in [kWh/$^\circ$C], and $ u $ denotes heat flux in [kW] from the electric heater. Solar irradiation $ \Phi $ [kW/m$^2$] enters the building interior through the effective window area $ A_{w} $ [m$^2$] and heats the building envelope that has an effective area $ A_{e} $ [m$^2$].  For the sake of notational simplicity, the subscript $i$ is not explicitly shown here. The parameters used to simulate the different houses are shown in Table~II.

\begin{figure}
	\centering
	\includegraphics[width=\linewidth]{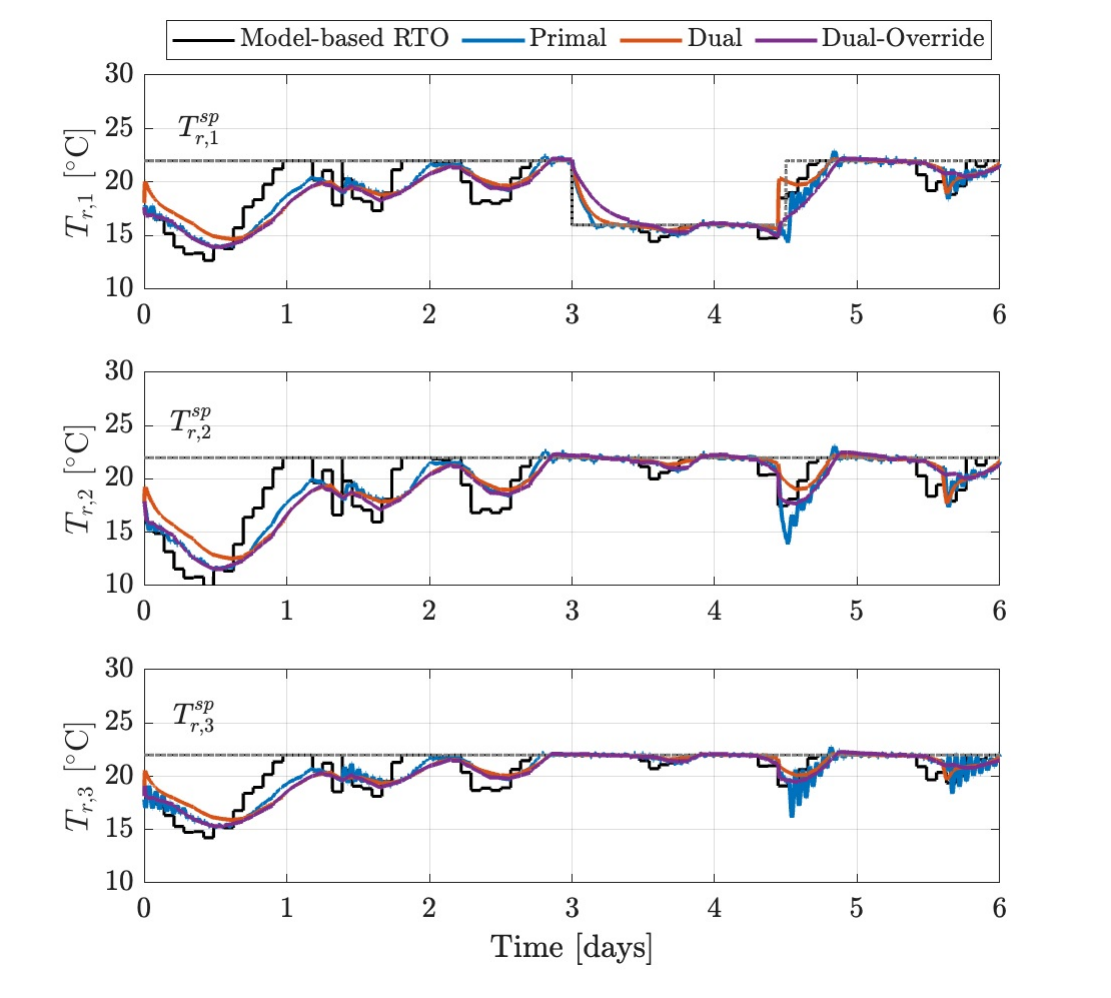}
	\caption{Case Study 1: Temperature profile for the three houses using primal-based (blue), dual-based (red) and dual-based with override (purple) control structure, compared with the model-based numerical optimization solution (in black).}\label{Fig:Temperatures}
\end{figure}
\begin{figure*}
	\centering
 \hspace*{-1.5cm}
	\includegraphics[width=1.12\linewidth]{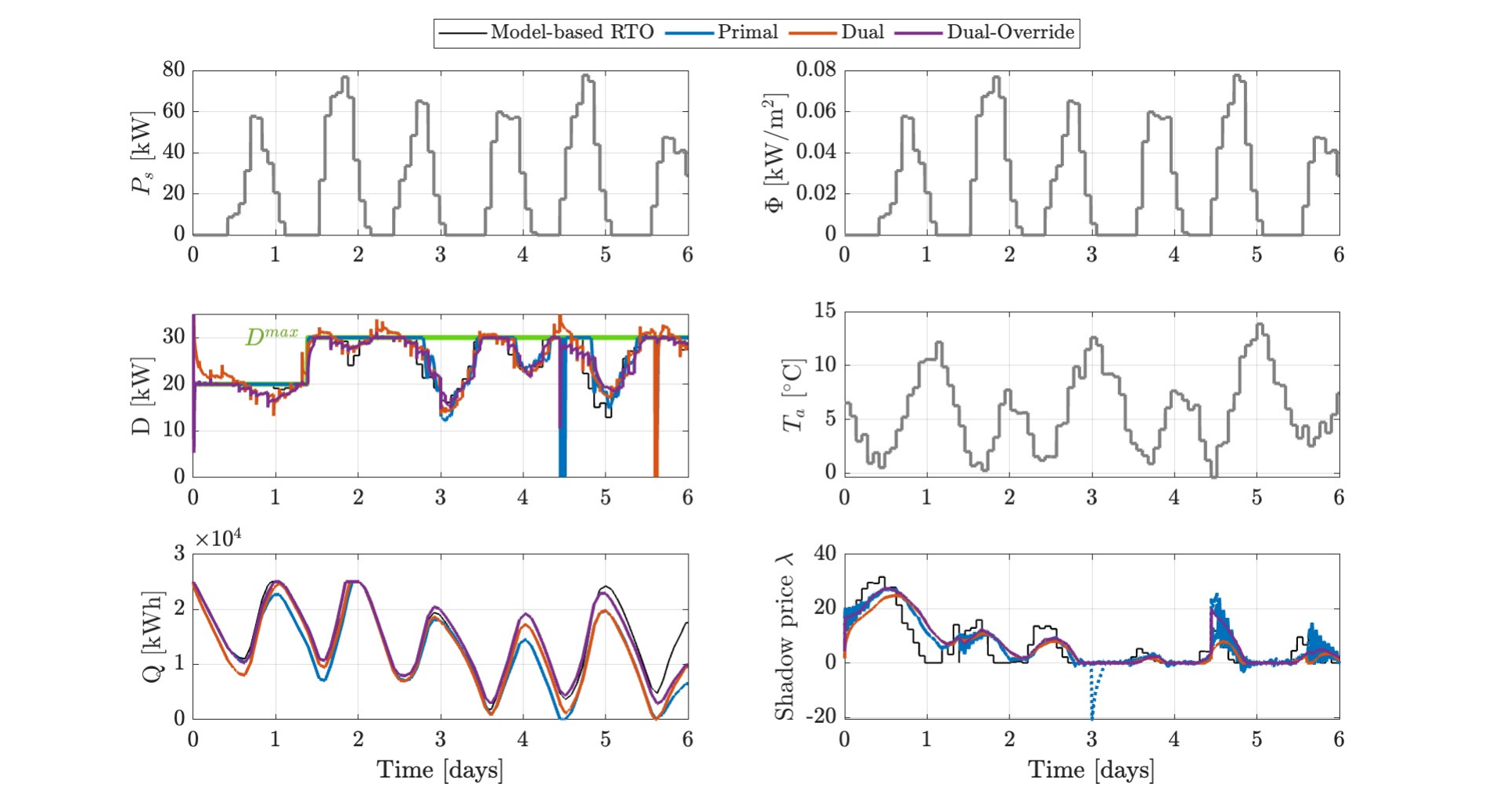}
	\caption{Case Study 1: Simulation results using primal-based (blue), dual-based (red) and dual-based with override (purple) control structure, compared with the model-based RTO solution (in black). Disturbances are shown in gray.}\label{Fig:Temperature2}
\end{figure*}
The state of charge of the battery storage $Q$ in [kWh] is described by $
 \dot{Q}=  P_{s} - D
$,
where  $D = \sum_{i=1}^{N}u_{i} $ is the total power discharged from the battery by the consumers, and $P_{s}= \omega \Phi$ is the amount of power generated from the solar panels to charge the battery, $\omega $ denoting the wattage of the solar panel and $\Phi$ denoting the solar irradiance.

In the context of the described small-scale energy hub with three residential buildings sharing power from a common solar energy source with battery storage for heating, demand response plays a crucial role in optimizing energy usage and ensuring grid stability. Each residential unit, modeled as an RC-thermal system, aims to maintain its interior temperature $T_{r}$ at a desired level $T_{r}^{sp}$ while minimizing deviations from this setpoint i.e., $J_{i}:= (T_{r,i} - T_{r,i}^{sp})^2$. This is achieved through the adjustment of heating settings and utilization of battery power to meet heating demands. However, to prevent excessive discharge and depletion of the battery, there are constraints on the total power that can be drawn from the battery across all buildings given by $\sum_{i=1}^{N}u_{i} \leq D^{max}$.

Achieving demand response automatically requires the local controllers within each residential unit to be designed such that it can adjust its electricity consumption to balance comfort requirements with energy availability and grid constraints, thereby optimizing overall system performance.
This section provides a comparative study of achieving demand response in this system  using  the three distributed feedback optimizing control architectures. The distributed feedback optimizing control architectures are also benchmarked using the model-based steady-state optimization approach, where each problem solves a numerical optimization problem online.

The three houses have a desired temperature of $T_r^{sp}$ = 22$^\circ$C. For 36 hours starting from day 3, house 1 has a reduced desired temperature of  $T_r^{sp}$ = 16$^\circ$C.
Given the ambient temperature $T_a$ and the solar irradiation $\Phi$ (shown in Fig.~\ref{Fig:Temperatures} in gray), solar power generation $P_s$  alone  is insufficient to maintain the desired temperature at all times.

Fig.~\ref{Fig:Temperature2} presents the results for battery discharge power (D) and state-of-charge (Q) under three control strategies: primal-based (blue), dual-based (red), and dual-based with override (purple). Fig.~\ref{Fig:Temperatures} shows the corresponding temperature profiles in each house.

The shadow price of the shared resource using dual-based (red), dual-based with override (purple), and the opportunity cost reported by the three houses (blue) are shown in the bottom right subplot in Fig.~\ref{Fig:Temperature2}. These are benchmarked against the model-based optimization solution (black), from which it can be seen that all the three methods are able to achieve asymptotic steady-state operation.

The differences between the three control structures are most evident when there is a change in the active constraints. To simulate this,  a step change in the value of $D^{max}$, the total power that can be drawn from the battery, is performed at some time after day 1. It should be noted that neither value is sufficient to drive the temperature to the required setpoint during night time (when the ambient temperature is low). Thus, the coupling constraint is still active at night, but become inactive during the day when the ambient temperature is higher.

Multiple instances can be observed where the total power drawn exceeds the constraint limit when using the dual-based control architecture. This occurs because the constraint enforcement operates on a slower timescale. When constraint violation is detected and the dual-override controller is activated, the fast controller intervenes by reducing the value of \textbf{u}$_{1}^g$. This results in an increase in the shadow price $\lambda$ on a slower time scale. This phenomenon is depicted in \textcolor{black}{Fig.~\ref{Fig:Temperature2}}, where the value of the dual variable (shadow price) is higher for the dual-override method than the dual-only method when the constraint is active. By design, the primal method (blue) ensures constraint satisfaction at all times. However, tuning the controllers for similar closed-loop response times led to overly aggressive control actions.

Interestingly, when the constraint is active, the shadow price of the dual-override method (purple) aligns with the primal-based method (blue) as shown in \textcolor{black}{Fig.~\ref{Fig:Temperature2}}. However, when the constraint is not fully utilized, the shadow price of the dual-override method (purple) resembles closely the dual-based method (red). The behavior of the shadow price in the dual-override method seems to be a hybrid of the primal and dual-based approaches. It prioritizes constraint satisfaction (like the primal method) when the constraint is active but switches to a more cost-oriented behavior (like the dual-based method) when there's spare capacity. 

\section{Case Study II - Optimal resource allocation in a subsea production network}\label{Sec: Lab Setup}

\subsection{Problem setup}
Subsea production systems involve extracting hydrocarbons from underground reservoirs through seabed wells and transporting them via pipelines to the surface. In cases of low reservoir pressure, methods like gas-lift are used, injecting compressed gases into the well tubing to reduce fluid density and minimize pressure losses. However, excessive gas injection can increase frictional pressure drop, necessitating optimal allocation of lift gas among wells to maximize production.
\begin{figure*}
	\centering
	\includegraphics[width = 0.9\textwidth]{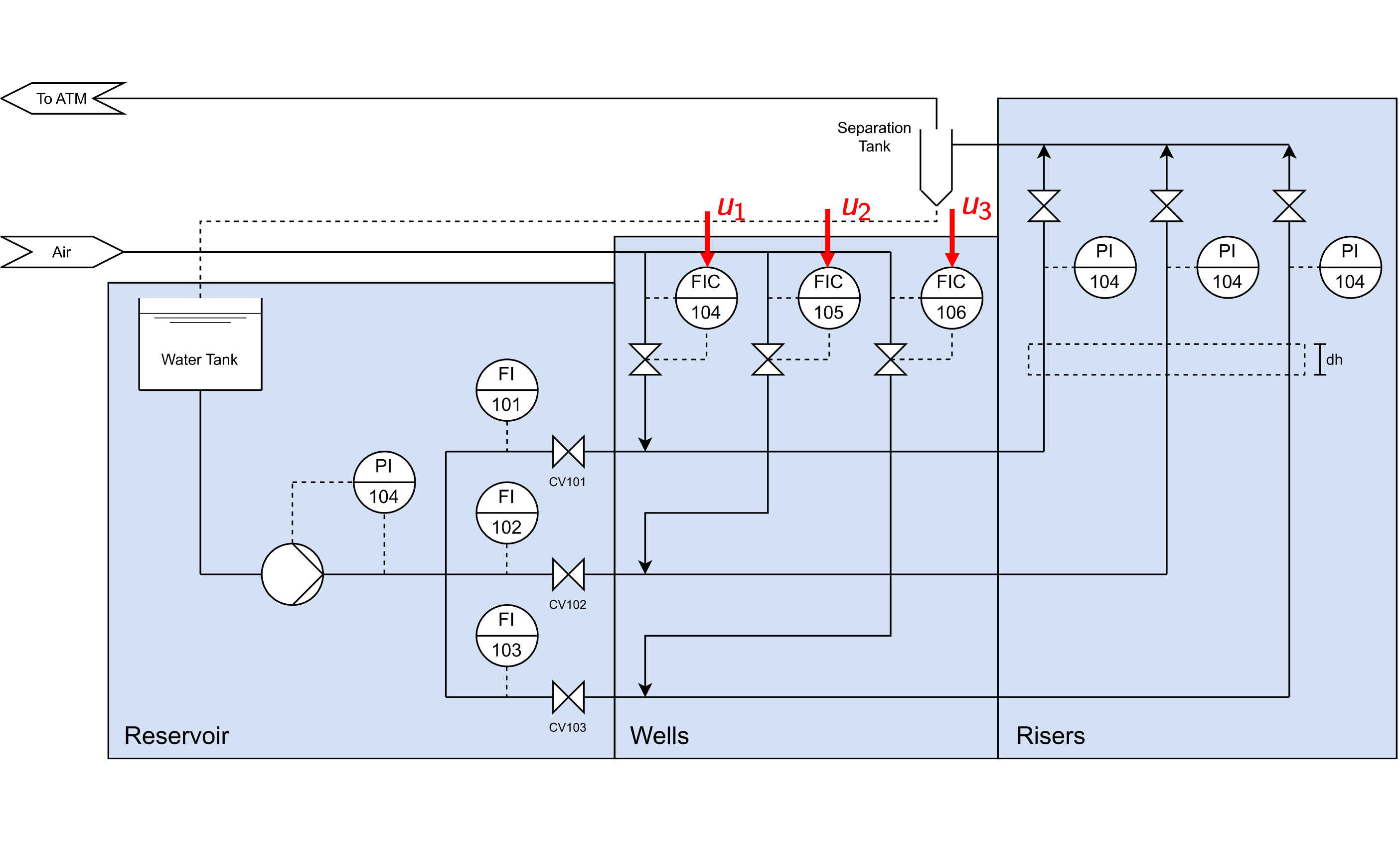}
	\caption{Case Study 2: This experiment's schematic is based on \cite{Matias2022Steady-stateRig}. The system measurements, denoted as $\mathbf{y}_\mathbf{p}$, include the well top pressures (PI101, PI102, and PI103), the pump outlet pressure (PI104), the liquid flowrates (FI101, FI102, and FI103), and the gas flowrates (FI104, FI105, and FI106). To control the gas flowrates, denoted as $\mathbf{u}=\begin{bmatrix} {Q}_{gl,1}& {Q}_{gl,2} & {Q}_{gl,3} \end{bmatrix}^{\top}$ to the calculated setpoints denoted as $\mathbf{u}^{sp}=\begin{bmatrix} {Q}_{gl,1}^{sp}& {Q}_{gl,2}^{sp} & {Q}_{gl,3}^{sp} \end{bmatrix}^{\top}$, three PI controllers are used. The reservoir valve openings (CV101, CV102, and CV103) represent system disturbances and vary during the experiments to simulate different reservoir behaviors. A PI controller maintains a constant pump outlet pressure throughout the experiment.}
	\label{fig:Lab_rig_schematic}
\end{figure*}
To emulate a subsea gas-lifted oil production system, we utilize a laboratory-scale experimental rig housed in the Department of Chemical Engineering at NTNU. This rig operates using water and air as working fluids instead of oil and gas for simplicity, without compromising the observation of the gas lift phenomenon. Hence, it is well-suited for studying production optimization techniques focused on gas lift effects. Illustrated in Fig. \ref{fig:Lab_rig_schematic}, the system comprises a reservoir, well, and riser section.
\begin{itemize}
\item  In the reservoir section, a stainless steel tank, centrifugal pump, and three control valves (CV101, CV102, and CV103) mimic reservoir disturbances like pressure oscillations or depletion. The reservoir produces only liquid, with outflow rates ranging from 2 L/min to 15 L/min, monitored by flow meters (FI101, FI102, and FI103). A PI controller maintains the pump's outlet pressure (PI104) at a constant 0.3 barg.
\item The well section includes three parallel flexible hoses (2 cm inner diameters, 1.5 m length) where pressurized air (approximately 0.5 barg) is injected by three air flow controllers (FIC104, FIC105, and FIC106) within a range of 1 sL/min to 5 sL/min, approximately 10 cm after the reservoir valves.
\item The riser section comprises three vertical pipelines (2 cm inner diameters, 2.2 m height), orthogonal to the well section, with pressure measured at the top (PI101, PI102, and PI103). After sensors, three manual valves are kept open during experiments to vent air to the atmosphere and recirculate liquid to the reservoir tank. Additional details are available in \cite{Matias2022Steady-stateRig}.
\end{itemize}

In this experimental setup, the optimization problem aims to maximize the network liquid flow rate, which is the combined production of liquid from three wells, while taking into account a limited amount of gas-lift available for injection. To express the economic objectives in line with problem \eqref{Eq:gen_prob}, we can state the following:

\begin{equation}    \label{Eq:rig_obj}
	\begin{split}
		&\sum_{i=1}^3 J_i(\mathbf{u}_i,\mathbf{d}_i)=\\
		& -20 Q_{l,1}\left(\mathbf{u}_1,\mathbf{d}_1\right) - 25 Q_{l,2}\left(\mathbf{u}_2,\mathbf{d}_2\right) - 30 Q_{l,3}\left(\mathbf{u}_3,\mathbf{d}_3\right)
	\end{split}
\end{equation}
where $Q_{l,i}$ denotes the produced liquid flowrate from well $i$. For illustration, we assume that the wells have different hydrocarbon prices as shown above. The production from the wells are controlled using the injected gas-lift flowrates denoted by $\mathbf{u}_{i} = {Q}_{gl,i}$, which are the decision variables for the optimization problem.  These are provided as setpoints to flow controllers (FICs) $104$, $105$, and $106$ to maintain the air injection flowrates at their respective setpoints (shown in red in Fig.~\ref{fig:Lab_rig_schematic}).  
Additionally, the disturbances  $\mathbf{d}_i$ are emulated by adjusting the valve openings CV101, CV102, and CV103.

The total gas availability, which is a shared (input) constraint, can also be expressed as follows,
\begin{equation}\label{lab_const}
	\begin{split}
	 \sum_{i=1}^3 \mathbf{g}_i(\mathbf{u}_i,\mathbf{d}_i) - \bar{\mathbf{g}} =  \sum_{i=1}^3 Q_{gl,i}  - Q_{gl}^{max}
	\end{split}
\end{equation}
\begin{remark}\label{Remark 4}
	Note that the total gas availability in general is an inequality constraint. However, in this experimental setup, we found that the gas lift constraint is always active. Hence the selector block does not select the marginal controller in the last subsystem in our simulation results.
\end{remark}
We now experimentally validate the three distributed feedback-optimizing control structures described in Section~\ref{Sec: Problem Formulation} for optimal resource allocation in this subsea production network, where each well represents a subsystem. 
Each subsystem employs an extended Kalman filter for local gradient estimation.  The local controllers are tuned using SIMC-rules \cite{Skogestad2003SimpleTuning}, while ensuring that the controllers in the slow timescale are 5-10 times slower than the controllers in the fast timescale (cf. Table~IV).  Prior to the implementation in the rig, we validated the controller tunings using a high-fidelity dynamic MATLAB model that includes lower-layer controller dynamics, i.e. FICs, and tuned noise parameters based on rig data. The  controller and tuning parameters used in this work are shown in Table  \ref{tab:tuning_main}. Detailed parameters are available on our Github page\footnote{{https://github.com/Process-Optimization-and-Control/ProductionOptRig}\label{foot}}, and a comprehensive model description can be found in \cite{Matias2022Steady-stateRig}. 	The simulator model is solely utilized for controller tuning parameter validation, preceding its implementation on the physical rig. All results in Section~5 originate from experiments conducted on the actual rig, not the simulator.

\subsection{Experimental Results}\label{Sec: Experiment}
We implemented the control structure design using dual decomposition, dual decomposition with override, and primal decomposition. When using the control structure design using primal decomposition, the  individual subsystems do not require local PI controllers as we consider shared input constraints (see also Remark~\ref{rem:primal}. Coordinator/central controllers in the slow timescale, responsible for equalizing all Lagrange multipliers, generate $\mathbf{t}_{i} = Q_{gl,i}^{sp}$, which are used by the flow controller FIC$_i$ to regulate air injection valve openings. 

\begin{figure}
	\centering
	\includegraphics[width=\columnwidth]{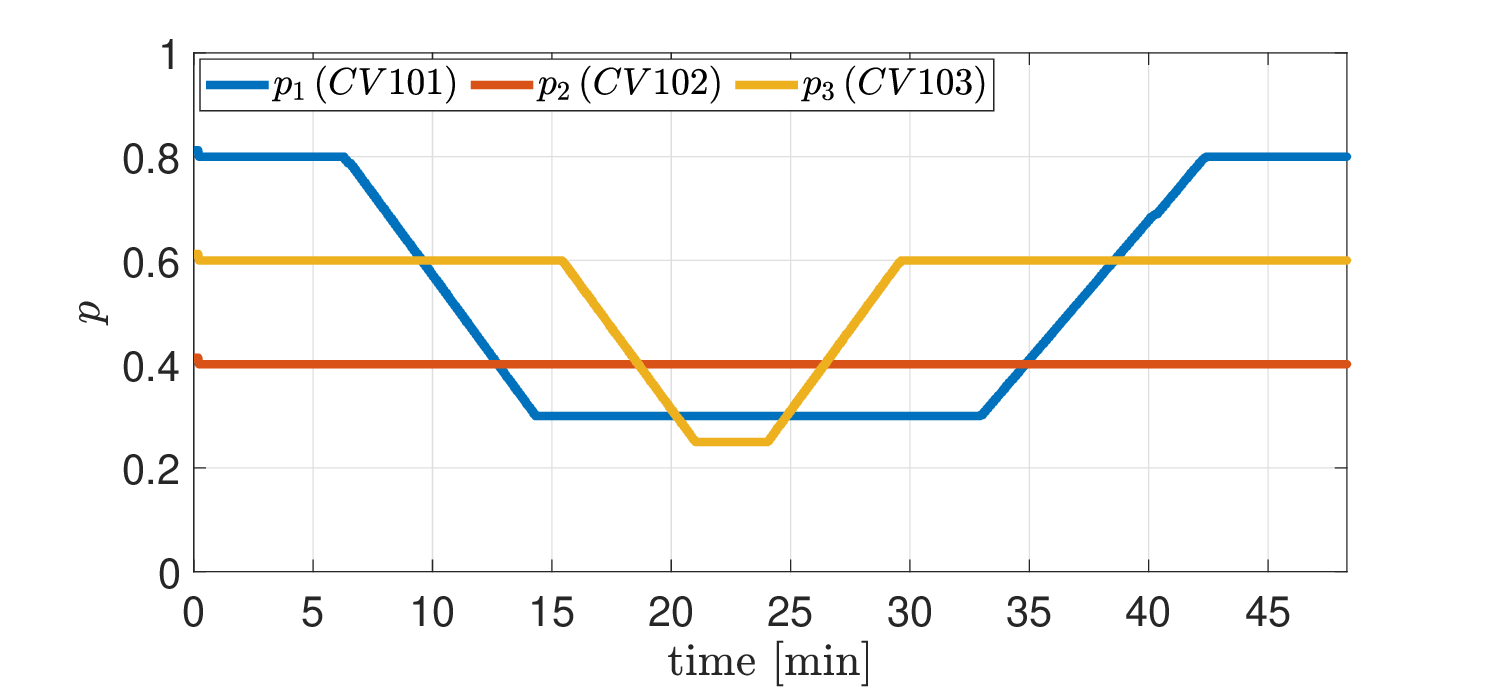}
	\caption{Case Study 2: The change of reservoir valve openings (CV101, CV102, and CV103) during the experiments for representing different reservoir behaviors. These reservoir valve openings are system disturbances.}
	\label{fig:Disturbance_Experiment}
\end{figure}

\begin{figure}
	\centering
	\includegraphics[width=\columnwidth]{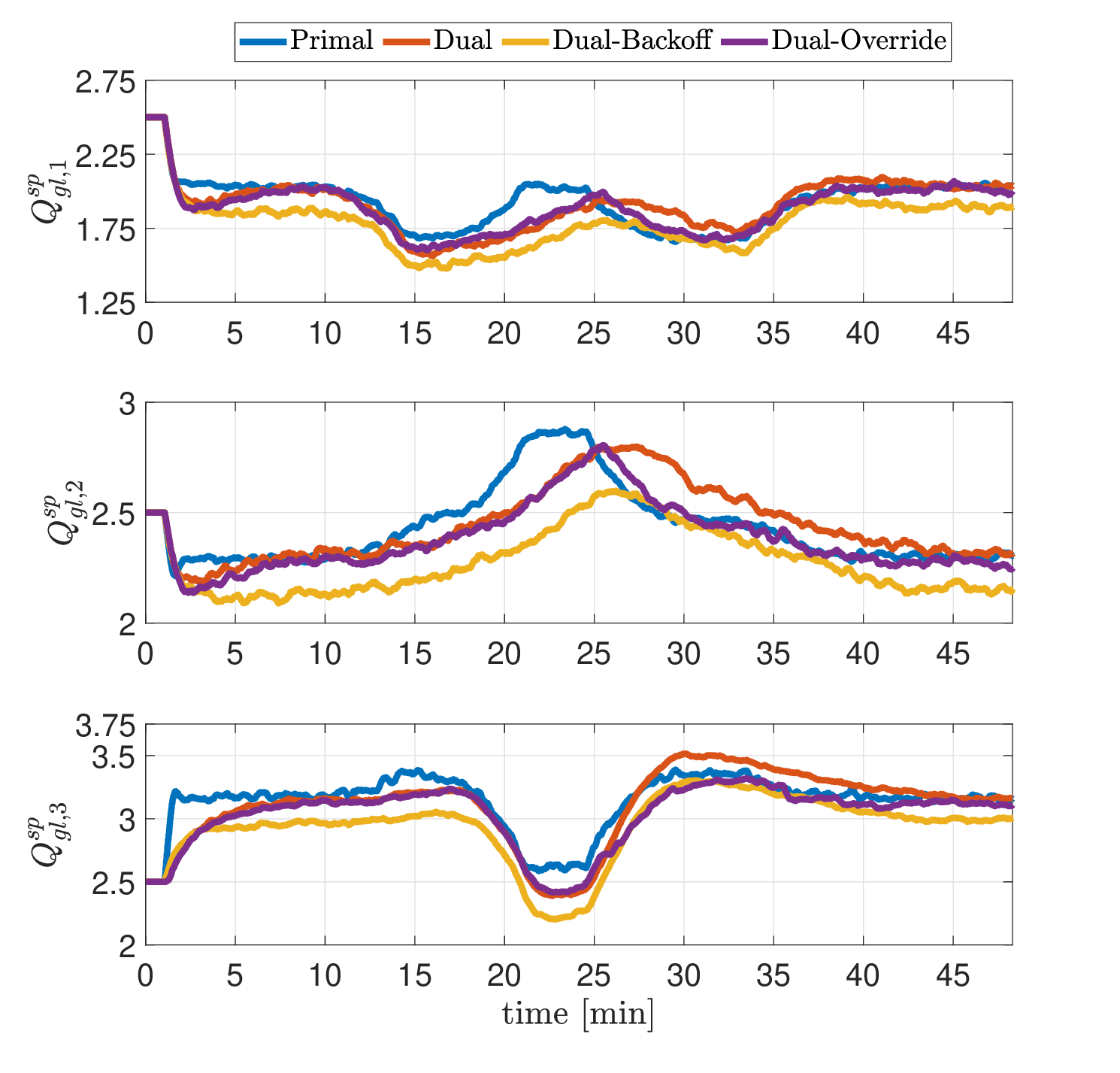}
	\caption{Case Study 2: The gas-lift flow rate setpoint $\left(\mathbf{u}^{sp}=\mathbf{Q}^{sp}_{gl}\right)$ of every wells due to reservoir parameter changing (disturbance) from the Experimental Lab Rig.}
	\label{fig:Input_SP_Comparative_Experiment}
\end{figure}

\begin{figure}
	\centering
	\includegraphics[width=\columnwidth]{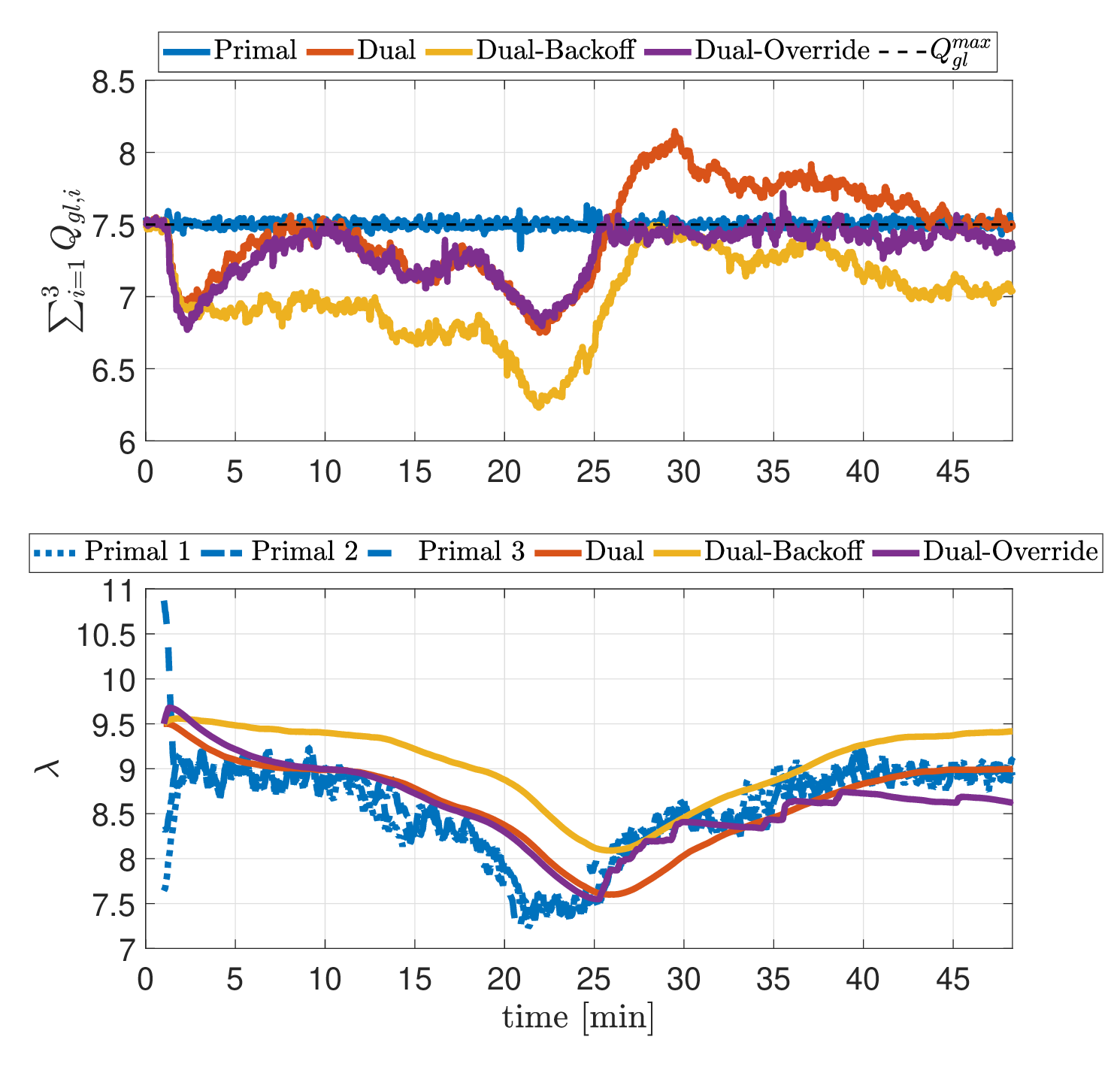}
	\caption{Case Study 2: Constraint satisfaction and Lagrange multiplier evolution due to reservoir parameter changing (disturbance) from the Experimental Lab Rig.}
	\label{fig:SP_Comparative_Experiment}
\end{figure}

\begin{figure}
	\centering
	\includegraphics[width=0.85\columnwidth]{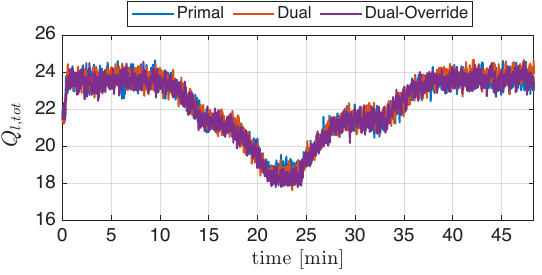}
	\caption{Case Study 2: Total produced liquid rate.}
	\label{fig:liquid rate}
\end{figure}

Fig. \ref{fig:Disturbance_Experiment} shows the reservoir valve openings ($CV101$, $CV102$, $CV103$) that we consider as the disturbance $\mathbf{d}_i$ in this experiment. The first disturbance occurs when the opening of $CV101$ gradually decreases from $t=6.5$ to $t=14$ minutes. We expect a decrease in the gas-lift injection in well 1, and a redirection of the gas supply to the other wells. The second disturbance occurs when the opening of $CV103$ also gradually decreases from $t=15.5$ to $t=21$ minutes. We expect that the gas supply to well 3 reduces with larger rate since the "hydrocarbon price" of this well is higher. Meanwhile the other wells will obtain more gas supply with larger rate as well. The third disturbance occurs when the opening of $CV103$ gradually increases from $t=24$ to $t=29.5$ minutes. We expect a reverse reaction compared to the second disturbance. Finally, the fourth disturbance occurs when the opening of $CV101$ also gradually increases from $t=33$ to $t=42$ minutes. Similarly, we expect a reverse reaction compared to the the first disturbance. We try to avoid sudden disturbance to ensure that the controller can adjust the plant smoothly. In the rig, we used a programming environment (LABVIEW \cite{Bitter2006LabVIEW:Techniques}) to automate the implementation of these disturbance. Therefore, it is possible to repeat the independent experiments with the same disturbance profile. The following results are the average of three independent experiments.
The experimental results comparing the three control structures are presented in Fig. \ref{fig:Input_SP_Comparative_Experiment}-Fig. \ref{fig:liquid rate}, where the results from the dual-based scheme are shown in red, dual with backoff is shown in yellow, dual with constraint override is shown in purple, and primal-based scheme is shown in blue. 


\textit{Comparison of constraint satisfaction:} As seen in the top plot of Fig. \ref{fig:SP_Comparative_Experiment}, controlling the coupling constraint in the slow timescale  in the dual-based scheme, leads to  significant violation of the coupling constraint. If this constraint violation is not desirable, one would have to impose a back-off for the constraint control. However adding a back-off can lead to steady-state performance loss. To address this, we proposed the dual-based scheme with constraint override. In Fig. \ref{fig:SP_Comparative_Experiment}, we can see that when the constraint is not active, this method (shown in purple) gives the same response as the dual-based approach (in red). However, at $t =25$, when the coupling constraint is violated, the constraint controller in the fast timescale overrides the gradient controller, and hence we see that the constraints are satisfied in the fast timescale.  In contrast,  the control structure based on primal decomposition ensures feasibility of the coupling constraints at all times as expected.  

The bottom plot of Fig. \ref{fig:SP_Comparative_Experiment} shows the associated Lagrange multiplier for the different control structures. It is important to notice that all the three local Lagrange multipliers of the primal-based approach converge to the same trajectory, indicating that the marginal cost is equal for all the subsystems. Meanwhile, each dual-based scheme has one Lagrange multiplier, and they are slower as they are controlled in the slow timescale. It is interesting to notice that before $t=25$ minute, dual-based with override has similar performance (trajectory) to dual-based scheme. However, after $t=25$ minute, dual-based with override has quite similar performance (trajectory) to primal-based scheme. Note that this phenomenon was also observed in case study 1 (cf. Section~\ref{sec:Demand response}). Specifically at the period between $t=25$ to $t=35$ minute, we can notice that dual-based with override is "enforcing" the trajectory (both the constraint and the associated dual variable) in order to minimize the constraint violation (in terms of magnitude and/or duration). The toal produced liquid rate using the different approaches are shown in Fig.~\ref{fig:liquid rate}.

\begin{figure}
	\centering
	\includegraphics[width=\columnwidth]{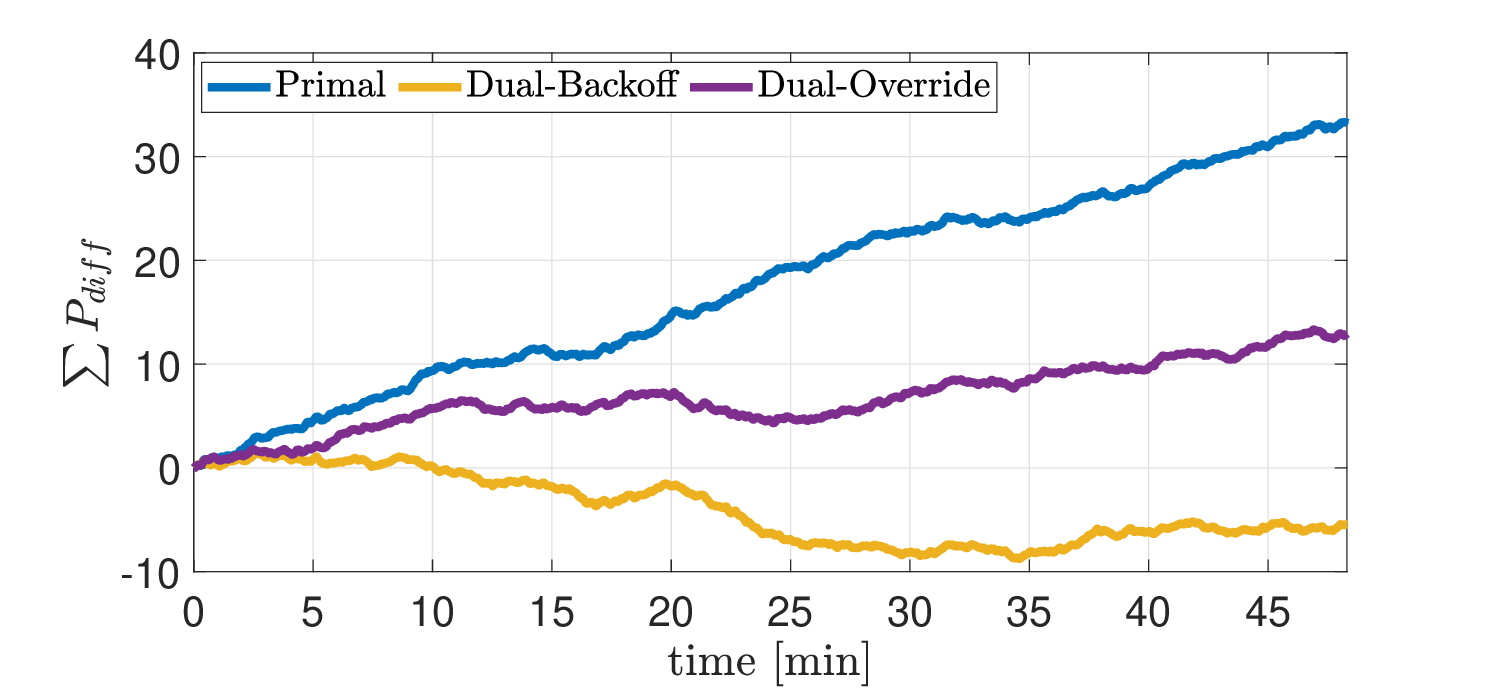}
	\caption{Case Study 2: The cumulative profit difference from the Experimental Lab Rig.}
	\label{fig:Prof_ave}
\end{figure}
To analyze the optimization performance of the three control structures, we compare the profit  obtained by the three schemes with the naive approach, where we consider fixed inputs, i.e., $\mathbf{u} = \begin{bmatrix} {Q}^{sp}_{gl,1}& {Q}^{sp}_{gl,2} & {Q}^{sp}_{gl,3} \end{bmatrix}^{\top} = \begin{bmatrix} \frac{{Q}_{gl}^{max}}{3} & \frac{{Q}_{gl}^{max}}{3} & \frac{{Q}_{gl}^{max}}{3} \end{bmatrix}^{\top}$. The naive approach illustrates the case in which no information about the system is available. Hence, the best alternative is to divide the available gas equally among the wells. This approach consistently satisfies the shared input constraint \eqref{lab_const}. It serves as a benchmark to demonstrate that the optimization schemes are capable of generating higher profits compared to not implementing any optimization at all.

Note that we do not include dual-based in comparison as it significantly violates the constraint, which implies that it has significant trajectory in infeasible region. In other words, dual-based requires significant extra resource starting from $t=25$ minutes (See Fig.~\ref{fig:SP_Comparative_Experiment}). To measure the performance, we plot the difference, in percentage, between the instantaneous profit of the approach of interest and the naive approach. The difference is calculated as
$
	P_{diff} = \frac{P-P_{naive}}{P_{naive}}.100
$
where $P$ is the profit of the approach of interest, and $P_{naive}$ is the profit of the naive approach.

Fig. \ref{fig:Prof_ave} depicts the cumulative profit/loss difference $\left(\sum P_{diff}\right)$ of all the three schemes of interest (i.e., primal-based, dual-based with back-off, and dual-based with override) to the naive strategy. Fig. \ref{fig:Prof_ave} shows that the primal-based scheme has the highest cumulative profit difference, and followed by dual-based with override scheme. This phenomenon can be attributed to the fact that the primal-based scheme consistently utilizes all available resources throughout the duration of the experiment. Conversely, the dual-based scheme with override takes more time, at the first 25 minutes, to fully utilize all available resources. Interestingly,  the dual-based scheme with back-off performed worse than the naive approach. 


\section{Discussion}\label{sec: Discussion}

\subsection{Controller placement and information exchange}
A key advantage of the proposed distributed feedback optimizing control architectures lies in their inherent modularity.  In addition to the horizontal decomposition across the different subsystems, the feedback optimizing control architectures presented in this work  utilize a two-layer, cascaded structure (vertical decomposition). In the faster layer, local controllers optimize performance within their subsystems based on the provided coordinating information. In the slower layer, the coordinating information is updated and  communicated to the lower layer  to achieve optimal operation for the entire system.

The nature of the coordinating information differs between control structures. Dual-based and dual-based with override architectures use shadow prices as the coordinating information. These prices reflect the  cost of consuming the shared resource, allowing local controllers to dynamically adjust their behavior based on real-time conditions. The prices are adjusted by a central coordinator based on the reported consumption of the shared resource.  In contrast, the primal-based structure employs a resource-directive approach, where each subsystem receives a pre-allocated amount of the shared resource, which is adjusted based on the reported opportunity costs.

In all the feedback optimizing control structures discussed in this paper, local controllers in the fast timescale and gradient estimators are situated within each subsystem. This arrangement offers a significant advantage: different subsystems can independently execute local feedback optimizing controllers without the need to share all local models and real-time measurements. This setup allows for coordinated action among different subsystems with minimal information exchange. Notably, different subsystems can employ controllers with varying closed-loop time constants and sampling times, enhancing flexibility.

\subsection{Timescale separation}
{As seen in Fig.~\ref{Fig:Dual_CSD}-\ref{Fig:Primal_CSD}, the controllers in the upper layer send information to the controllers in the lower layer in a cascaded fashion. When designing control architectures in a csacade fashion, for the controllers in the fast timescale and the slow timescale to not interfere with each other, we need a clear timescale separation between the two layers.} This time scale separation is required for all the three control architectures.  This can be achieved by tuning the controller gains such that the closed-loop time constant of the lower level controllers are at least 5 times faster than the closed-loop time constant of the upper level controllers. {A detailed analysis of the timescale separation when designing feedback optimizing control archutectures can be seen in  \cite[Section 2.5]{Skogestad2023AdvancedElements}.}

In price-based coordination schemes (dual-based and dual-based with override), a central coordinating controller updates shadow prices in a slower timescale, shared by all subsystems. This coordinating controller broadcasts the updated shadow prices to each subsystem. In order to satisfy the timescale separation requirement, the closed-loop time constant of the coordinating controller that updates the shadow price is determined by the closed loop time constant of the slowest subsystem.

Conversely, in the primal-based control architecture, the marginal cost controller responsible for updating the allocation of the shared resource $\mathbf{t}_{i}$ in the slow timescale is also positioned within each subsystem. Consequently, the closed-loop time constants of marginal cost controllers that works in the slow timescale can be independently chosen within each subsystem, provided they are at least 5 times slower than the closed-loop time constants of their corresponding lower-level controllers. This arrangement is illustrated in Fig.~\ref{Fig:Primal_CSD}, where controllers in the slow timescale and fast timescale are located within the information boundary of each subsystem (indicated by a grey box).

\begin{remark}\label{rem:primal}
This work also reveals an interesting connection between distributed feedback optimizing control based on primal decomposition and existing control structures for parallel operating units. When the coupling constraint solely affects the input (i.e., $\mathbf{g}_i(\mathbf{u}_i,\mathbf{d}_i) = \mathbf{u}_i$
 ), the primal decomposition framework simplifies the control architecture to solely utilize marginal cost controllers within each subsystem. The resulting control structure is identical to the  control structure employed for optimizing parallel operating units, as demonstrated in \cite{downs2011industrial,DK2020ESCAPE,DK2019CEP}. This finding highlights the general applicability of the proposed distributed feedback optimizing control using  primal decomposition framework, as it encompasses previously known control structures within the process control domain.
\end{remark}

{The general approach of using vertical decomposition with time-scale separation used in this paper can also be found in other works in the literature, such as in two-layer hierarchical MPC, which also deals with the dynamics \cite{nigro2024distributed,la2023predictive}, where the lower layer controllers are implemented at fast sampling rates, and the upper layer controllers are implemented in slow sampling rates.}

\subsection{Pairing}
When designing feedback controllers, an important consideration is how to pair the controlled variables and the manipulated variables. For dual-based approach, the pairing is rather straightforward, where for the local controllers, it is natural to pair diagonally, and for the coordinator in the slower timescale, we pair the constraint to its corresponding Lagrange multiplier (shadow price).

In the dual-based approach with constraint override, we adhere to the same pairing structure as the dual-based control. However, we must also identify the critical subsystems where the local gradient controllers will be overridden by the constraint controllers, and determine the pairing of the constraint controllers. It's essential to avoid pairing the constraints with manipulated variables that quickly saturate, as this would compromise our ability to control the constraint in the fast timescale. This can be done by pairing the constraint with the inputs with a large process gain.

For the primal-based control architecture, one must choose the controller pairings for the local constraint control and gradient controllers. This may not always be trivial since there is no systematic approach that provides a unique pairing. In general, the  pairing must be chosen such that the steady- state relative gain array (RGA) of the resulting transfer matrix is non-negative and close to identity matrix at crossover
frequencies \cite{Skogestad2005MultivariableControl}. {One could also use the Niederinski index (NI) for a high-level stability check of the proposed pairing (e.g., NI$<$0 indicates likely instability for the chosen pairing). However RGA  provides more detailed information about loop interactions, helping refine the pairing and providing insight into which variables will strongly influence each other.} Furthermore, to avoid the constraint controllers from saturating, one must pair the constraints with inputs with a large process gain. In practice, if the inputs of the constraint controller gets saturated, additional strategies, such as split-range control, must be employed to effectively implement the primal-based scheme. In order to minimize time delay, it is also desirable to pair manipulated variables that are physically close to the controlled variables (pair-close rule). {While Grammians are not typically used for direct pairing decisions, they can be valuable when designing control systems that must ensure strong controllability and observability, making them more applicable in situations where the dynamics of the system are critical, and are generally more useful in state-space-based control designs.}

\subsection{Coupling constraints and local constraints} According to Assumption~\ref{asm:sufficientMVs}, the primal-based approach is capable of handling a maximum number of coupling constraints equal to or less than the minimum number of local manipulated variables $\left(m\leq  \min_i n_{i}\right)$. In contrast, the dual-based approach, is not limited by the number of coupling constraints. This fundamental distinction implies that the dual-based scheme has the advantage of being able to accommodate a larger number of coupling constraints compared to the primal-based method.

In addition to the  coupling constraints, the subsystems may also have additional local constraints $\mathbf{h}_{i}(\mathbf{u}_{i},\mathbf{d}_{i}) \leq 0$, with $\mathbf{h}_{i} \in \mathbb{R}^{m_{i}}$.  The presence of local constraints can complicate the primal-based control structure design.  Firstly, the primal-based approach necessitates an adequate number of manipulated variables such that $n_{i}\geq m+m_{i}$, ensuring that the available manipulated variables can manage both the local constraints and the coupling constraint $\mathbf{g}_{i} - \mathbf{t}_{i}$. Moreover, the local constraints may conflict with the coupling constraint for a given $\mathbf{t}_{i}$. Consequently, the subsystems would need to share information regarding the local constraints $\mathbf{h}_{i}$ with all other subsystems to ensure that the allocated $\{\mathbf{t}_{i}\}$ remain consistent with the local constraints, which may raise concerns regarding data privacy.

 In contrast, the dual-based  control structures (with or without constraint override), does not suffer from the limitations imposed by additional local constraints. The local constraints can be treated in the same form as coupling constraints,  in the sense that they are relaxed and controlled in the slower timescale using the  corresponding Lagrange multiplier. In this case, in addition to the common coordinating controller that control the coupling constraint, the slow timescale also includes the local  constraint controllers within each subsystem that updates the corresponding Lagrange multiplier $\mu_i (t+1) = \max\left(0,\mu_i(t) + K_{h_{i}}\mathbf{h}_{i}\right)$. In the fast timescale, the local controllers control the self optimizing variable  $ \mathbf{c}_{i}(\lambda,\mu_i) = \nabla J_{i} + \nabla\mathbf{g}_{i}\lambda + \nabla \mathbf{h}_{i}\mu_i$.  This characteristic distinguishes the dual-based approach, which proves to be more versatile and adaptable in the presence of various types of constraints within a distributed feedback-optimizing control architecture.

 To summarize, the primal-based control architecture suffers from poor scalability in terms of the number of constraints that can be handled, whereas the dual-based control architectures can be easily scaled to handle any number of local and coupling constraints.
\subsection{Gradient estimation}
The feedback optimizing control structures explored in this study necessitate online estimation of steady-state cost and constraint gradients around the current operating point. Here, we employed a model-based gradient estimation approach using an extended Kalman filter (details of which can be found in \cite[Appendix A]{Dirza2022ExperimentalRig}. However, it's worth noting that any model-based or model-free gradient estimation method could potentially be utilized with the  control architectures described in this paper. We refer the reader to \cite{Krishnamoorthy2022Real-TimeReview} for a comprehensive overview of various model-based and model-free gradient estimation methods that can be used for feedback optimizing control.

While model-free gradient estimation schemes may appear appealing to address plant-model mismatch, caution is advised due to the timescales involved with estimating the steady-state gradients. For instance, sinusoidal perturbation-based gradient estimation (commonly employed in extremum-seeking control approaches) requires clear timescale separation between plant dynamics, sinusoidal perturbations, and local controller convergence, which can  be prohibitively slow, especially for systems with slow dynamics \cite{Krishnamoorthy2022Real-TimeReview,Srinivasan2019i110thSchemes,nesic2009extremum}. In addition, achieving the timescale separation necessary for coordination  may further impede overall system convergence. 

In the dual-based control structures (with and without constraint override), plant-model mismatch can be handled by using a combination of model-based gradient estimation (fast timescale) with model-free gradient estimation (slow timescale), where the local controllers control the self-optimizing variable $\mathbf{c}_{i}(\lambda)$ to a setpoint of $\delta_i := \nabla_{p}\mathcal{L} - \nabla\mathcal{L}$ instead of zero. This essentially corrects for the plant-model mismatch in the slow timescale by modifying the setpoint of the local gradient controllers. This is akin to the feedback variant of modifier adaptation scheme for steady-state real-time optimization \cite{marchetti2009modifier,milosavljevic2016optimal}.

On the other hand, the use of such modifier terms for handling plant-model mismatch are not trivial in the primal-based control structure. The primal-based approach estimates local Lagrange multipliers by evaluating a matrix calculation as depicted in equation \eqref{Eq: local lang}. This is done using the estimated gradients, and is therefore susceptible to plant-model mismatch. Furthermore, this approach relies on the  assumption that the local Lagrange multipliers are both unique and exist. In contrast, the dual-based approach, whether with or without an override, estimates the Lagrange multipliers through feedback mechanisms using the measured constraints, thereby circumventing the need for such matrix calculations. More importantly, the Lagrange multipliers computed by the coordinating controllers are not susceptible to plant-model mismatch. As a result, the dual-based method exhibits a reduced likelihood of encountering numerical issues, offering a more robust alternative in more general practical applications.

\subsection{Extension to other types of network interconnections }
Unlike the price-based coordination mechanism, which requires a central entity to update the shadow prices,  the primal-based control architecture  naturally extends towards a  peer-to-peer coordination  framework, where instead of controlling $\lambda_N - \lambda_i$  for all $i=1,\dots,N-1$,  equal marginal cost can also be achieved by controlling $\lambda_j - \lambda_i$ where the pair $(i,j)$ indicates a communication edge between subsystems $i \in \{1,\dots,N-1\}$ and $j\in \{1,\dots,N\}$ in an arbitrary connected  graph network. Here the term \emph{connected graph} indicates that  there is a path between any subsystem $i$ and subsystem $N$.

\begin{figure}
    \centering
    \includegraphics[width=0.5\linewidth]{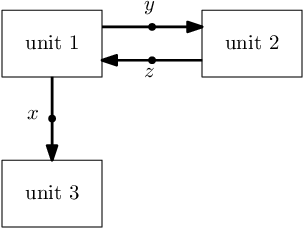}
    \caption{Schematic representation of a process system interconnected via flow streams.}
    \label{fig:connected}
\end{figure}

The three control system architectures presented in this paper can also be extended to  other problem classes with additively separable cost.
Consider for example the optimization problem with $N=3$  subsystems with interconnected flowstreams as shown in Fig.~\ref{fig:connected}, \[ \min_{\{\mathbf{u}_i\},x,y,z} \sum_{i=1}^NJ_i(\mathbf{u}_i,\mathbf{d}_i,x,y,z)\] where the common variables $x,y,z$ denotes flow stream between the different subsystems as shown in Fig~\ref{fig:connected}. This problem can be reformulated as
\begin{align*}
    \min_{\{\mathbf{u}_i,\mathbf{d}_i,x_i,y_i,z_i\}} & \;  \sum_{i=1}^N J_i(\mathbf{u}_i,\mathbf{d}_i,x_i,y_i,z_i) \\
    \text{s.t.} \; & x_1-x_3 =0 \\
    & y_1-y_2 = 0, \quad z_2-z_1 = 0
\end{align*} The control architectures detailed in this paper can then be applied to this case by setting $\mathbf{g}_1 = [x_1,y_1,-z_1]^{\mathsf{T}}$, $\mathbf{g}_2 = [0,-y_2,z_2]^{\mathsf{T }}$,$\mathbf{g}_3= [-x_3,0,0]^{\mathsf{T}}$ and $\bar{\mathbf{g}} = 0$, such that $ \sum \mathbf{g}_i - \bar{\mathbf{g}} = 0$ captures the flow stream mass balance constraints.


\begin{table*}[t]
\centering
\caption{Pros and Cons of Distributed Feedback Optimizing Control Architectures}
\begin{tabular}{|p{2cm}|p{7cm}|p{7cm}|}
\hline
\textbf{Control architecture} & \textbf{Pros} & \textbf{Cons} \\
\hline
Primal-Based Control & - Ensures constraint feasibility also during the transients. \newline- Naturally extends to a peer-to-peer coordination framework. \newline - Closed-loop time constants of the upper layer can also be chosen independently. & - Poor scalability due to the number of constraints that can be handled. \newline - Requires sharing local constraint information among subsystems, raising data privacy concerns. \newline - Complex pairing process without a systematic approach. \newline - Lagrange multiplers susceptible to plant-model mismatch and numerical issues. \\
\hline
Dual-Based Control & - Scalable and adaptable to local and coupling  constraints (Local constraints can be treated similarly to coupling constraints). \newline - Lagrange multipliers robust to plant-model mismatch. \newline - Pairing of controlled variables and manipulated variables is relatively straightforward.  & - Constraints handled at slow timescale. \newline  - Requires a central entity for updating shadow prices. \newline - Timescale separation needed between the central coordinating controller and local controllers dominated by the slowest subsystem. \\
\hline
Dual-Based with Override & - Inherits all advantages of the dual-based method. \newline - Capable of handling local constraints without additional complexity. \newline - Critical subsystems override local gradient controllers when necessary to handle constraints in fast timescale. & - Same central entity requirement as the dual-based method. \newline - Requires careful identification of critical subsystems for constraint override. \\
\hline
\end{tabular}
\newline
\begin{itemize}
\item All architectures leverage modularity through horizontal decomposition (across subsystems) and vertical decomposition (across timescales).
\item Timescale separation essential for avoiding interference between fast and slow layers in all architectures.
\item Minimal information exchange is required, enhancing system flexibility and independence.
\item Any model-based or model-free gradient estimation approach can be used with the control architectures.
\item  All methods can be extended to  problem classes with additively separable costs, ensuring versatility across different system configurations.
\end{itemize}

\label{table:pros_cons}
\end{table*}

\section{Conclusion}\label{Sec: Conclusion}

In this work, we validated three distributed feedback-optimizing control architectures through two distinct case studies. Our findings confirm that all three architectures can achieve optimal resource allocation at steady-state (i.e., asymptotic optimal operation). However, their transient behavior varies depending on the control structure and the tuning parameters of the controllers.
We have provided a comprehensive comparison of the pros and cons of each control architecture in Section~\ref{sec: Discussion}, summarized in Table~\ref{table:pros_cons}. These insights highlight the unique strengths and limitations of each approach.
Ultimately, the optimal choice of control architecture depends on the specific characteristics of the system, including the types and numbers of local and coupled constraints, the availability of manipulated variables, and the performance priorities. These findings offer valuable guidance for selecting the most suitable control strategy for diverse applications.


\appendices
\begin{table}[H]
	\caption{Building Parameters - Case Study 1}
	\begin{center}
		\begin{tabular}{c c c c c}
			\hline
			Description & Variable &  $i=1$ &  $i=2$ &  $i=3$ \\
			\hline
			\multicolumn{5}{c}{Heat Capacity [$^\circ$C/kW]} \\
			\hline
			Sensor & $C_{s}$ & 0.0549 & 0.0549 & 0.0549 \\
			Interior & $C_{i}$ & 0.0928 & 0.0835 & 0.1114 \\
			Envelope & $C_{e}$ & 3.32 & 2.656 & 4.648 \\
			Heater & $C_{h}$ & 0.889 & 0.889 & 0.889 \\

			\hline
			\multicolumn{5}{c}{Thermal Resistance between [$^\circ$C/kW]} \\
			\hline

			Interior \& Sensor & $R_{is}$ & 1.89 & 1.89 & 1.89 \\
			Interior \& Heater & $R_{ih}$ & 0.146 & 0.146 & 0.146 \\
			Interior \& Envelope & $R_{ie}$ & 0.897 & 0.7176 & 1.0764 \\
			Envelop \& Ambiance & $R_{ea}$ & 4.38 & 3.066 & 5.256 \\
			Interior \& Ambiance & $R_{ia}$ & 2.5 & 2 & 3 \\

			\hline
			\multicolumn{5}{c}{Area [m$^2$]} \\
			\hline

			Window Effective Area & $A_{w}$ & 5.75 & 3.8525 & 9.2 \\
			Envelope Effective Area & $A_{w}$ & 3.87 & 3.096 & 5.805 \\

		\end{tabular}
	\end{center}
	\label{tab:building params}
\end{table}

\begin{table}[H]
	\caption{Controller Parameters - Case Study 1} 
	\begin{center}
		\begin{tabular}{c c c}
			\hline
			Description & Variable & Value \\
			\hline
				Sampling time for simulation [s] & $t_{f}$ & 60 \\
				Sampling time for solar power [hr] & $t_{solar}$& 1.667 \\
			\hline
			\multicolumn{3}{c}{Primal Based} \\
			\hline
				Local Controller 1 I-Gain & $K_{I,1}$ & 0.052 \\
				Local Controller 2 I-Gain & $K_{I,2}$ & 0.052 \\
				Local Controller 3 I-Gain & $K_{I,3}$ & -0.1 \\
			\hline
			\multicolumn{3}{c}{Dual Based} \\
			\hline
				Coordinator I-Gain & $K_{I,slow}$ & 0.01\\
			 	Local Controller 1 P-Gain & $K_{P,1}$ & 3.85 \\
			 	Local Controller 2 P-Gain & $K_{P,2}$ & 2.5 \\
			 	Local Controller 3 P-Gain & $K_{P,3}$ & 2.5 \\
			 	Local Controller 1 I-Gain & $K_{I,1}$ & 0.01 \\
			 	Local Controller 2 I-Gain & $K_{I,2}$ & 0.01 \\
			 	Local Controller 3 I-Gain & $K_{I,3}$ & 0.01 \\

			\hline
			\multicolumn{3}{c}{Dual Based with Override} \\
			\hline
				Coordinator Anti Windup-Gain & $K_{AW, slow}$  & 0.01 \\
				Coordinator I-Gain & $K_{I, slow}$  & 0.1 \\
				Override Controller P-Gain & $K_{P, fast}$  & 0.1 \\
				Override Controller I-Gain & $K_{I, fast}$  & 0.03 \\
				Local Controller 1 P-Gain & $K_{P,1}$ & 2.5 \\
				Local Controller 2 P-Gain & $K_{P,2}$ & 2.5 \\
				Local Controller 3 P-Gain & $K_{P,3}$ & 2.5 \\
				Local Controller 1 I-Gain & $K_{I,1}$ & 0.01 \\
				Local Controller 2 I-Gain & $K_{I,2}$ & 0.01 \\
				Local Controller 3 I-Gain & $K_{I,3}$ & 0.01 \\
		\end{tabular}
	\end{center}
\end{table}

\begin{table}[H]
    \caption{Controller and Tuning parameters - Case study 2}
\begin{center}
\begin{tabular}{c c c}
 \hline
Description & Variable & Value \\ 
 \hline
 Exp. rig sensors sampling time & $T_{s}$ & {1} {second}\\
 \hline
 \multicolumn{3}{c}{Primal-based} \\
 \hline
Coordinator 1 I-Gain & $K_{I,1}$ & -0.016 \\ 
Coordinator 2 I-Gain & $K_{I,2}$ & -0.014 \\ 
Coordinator 3 I-Gain & $K_{I,3}$ & -0.011 \\ 
 \hline
 \multicolumn{3}{c}{Dual-based} \\
 \hline
        Coordinator I-Gain & $\alpha$ & {0.0088}                  \\ 
        Gradient Controller Input 1 I-Gain & $K_{I,1}$ &  0.016  \\ 
        Gradient Controller Input 2 I-Gain & $K_{I,2}$ & 0.014  \\ 
        Gradient Controller Input 3 I-Gain & $K_{I,3}$ & 0.011  \\ 
 \hline
 \multicolumn{3}{c}{Dual-based with Override} \\
 \hline
        Coordinator P-Gain & $K_{P,\alpha}$ & {1.7686}                 \\ 
        Coordinator I-gain & $\alpha$ & {0.00368}                 \\ 
        Gradient Controller Input 1 I-Gain & $K_{I,1}$ &  0.016 \\ 
        Gradient Controller Input 2 I-Gain & $K_{I,2}$ & 0.014  \\ 
        Gradient Controller Input 3 I-Gain & $K_{I,3}$ & 0.011  \\ 
        Override Controller I-Gain & $K_{I,3}$ & 0.1  \\         
        Override Controller AntiWindup-Gain & $K_{aw}$ & 0.4 \\
  \hline
 \multicolumn{3}{c}{Local set point controllers} \\
 \hline
        FIC 104 P-Gain & $K_{P,FIC,1}$ & 8560  \\
        FIC 105 P-Gain & $K_{P,FIC,2}$ & 8560  \\
        FIC 106 P-Gain & $K_{P,FIC,3}$ & 8560  \\
        FIC 104 I-Gain & $K_{I,FIC,1}$ & 100  \\
        FIC 105 I-Gain & $K_{I,FIC,2}$ & 100  \\
        FIC 106 I-Gain & $K_{I,FIC,3}$ & 50  \\
 \hline
\end{tabular}
\end{center}
\label{tab:tuning_main}
\end{table}

\bibliography{references}
\bibliographystyle{IEEEtran}

\vfill

\end{document}